\documentclass[1p]{elsarticle}
\usepackage{graphicx} 
\usepackage{subcaption}
\usepackage{mathtools}
\usepackage{amsmath}
\usepackage{amsfonts}
\usepackage{amssymb}
\usepackage{amsthm}
\usepackage{array}
\usepackage{xcolor}
\usepackage[colorlinks=true, allcolors=blue]{hyperref} 
\usepackage{cleveref}
\usepackage{footnote}
\usepackage[linesnumbered,ruled]{algorithm2e}
\usepackage{float}
\usepackage{verbatim}
\usepackage{diagbox}
\usepackage{array}
\usepackage{adjustbox}
\usepackage{geometry}

\crefformat{equation}{(#2#1#3)}
\Crefname{equation}{Equation}{Equations}
\Crefname{figure}{Figure}{Figures}
\Crefname{algocf}{Algorithm}{Algorithms}

\DeclareMathOperator{\R}{\mathbb{R}}
\DeclareMathOperator{\Risk}{\mathcal{R}}

\DeclareMathOperator{\E}{\mathbb{E}}

\DeclareMathOperator{\th_}{th}
\DeclareMathOperator{\diag}{diag}
\DeclareMathOperator{\tr}{tr}
\DeclareMathOperator{\MAP}{MAP}
\DeclareMathOperator{\LS}{LS}
\DeclareMathOperator{\DEIM}{DEIM}
\DeclareMathOperator{\prior}{prior}
\DeclareMathOperator{\post}{post}
\DeclareMathOperator{\noise}{noise}
\DeclareMathOperator{\opt}{opt}
\DeclareMathOperator{\logdet}{\log\det}
\DeclareMathOperator{\Dopt}{D-opt}
\DeclareMathOperator{\normal}{\mathcal{N}}
\DeclareMathOperator{\com}{\text{,}}
\DeclareMathOperator{\period}{\text{.}}
\DeclarePairedDelimiter\abs{\lvert}{\rvert}
\newtheorem{theorem}{Theorem}[section]

\newtheorem{lemma}[theorem]{Lemma}
\newtheorem{proposition}[theorem]{Proposition}
\newcommand{\id}{\mbox{d}}
\newcommand{\norm}[1]{\left \lVert #1 \right \rVert}
\newcommand{\sqbr}[1]{\left[#1\right]}
\newcommand{\Ord}[1]{\mathcal{O} \left(#1\right)}
\newcommand{\paren}[1]{\left(#1\right)}
\newcommand{\deter}[1]{\text{det}\left(#1\right)}

\newcommand{\Rank}[1]{\text{Rank}\left(#1\right)}

\newcommand{\bmat}[1]{\begin{bmatrix} #1 \end{bmatrix}}
\newtheoremstyle{definition}
  {10pt} 
  {10pt} 
  {} 
  {} 
  {\bfseries} 
  {.} 
  {5pt plus 1pt minus 1pt} 
  {} 
\theoremstyle{definition}
\newtheorem{definition}{Definition}[section]

\begin{document}

\begin{frontmatter}

\title{Bridging the Gap Between Deterministic and Probabilistic Approaches to State Estimation}

\author[a]{Lev Kakasenko} 
\author[a]{Alen Alexanderian}
\author[a]{Mohammad Farazmand\corref{cor1}}
\author[a]{Arvind K. Saibaba}
\cortext[cor1]{Email address for correspondence: \href{mailto:farazmand@ncsu.edu}{farazmand@ncsu.edu}}

\affiliation[a]{organization={Department of Mathematics, North Carolina State University},
            addressline={2311 Stinson Drive},
            city={Raleigh},
            state={North Carolina},
            postcode={27695},
            country={USA}}

\begin{abstract}
We consider the problem of state estimation from limited discrete and noisy measurements. 
In particular, we focus on modal state estimation, which approximates the unknown state of 
the system within a prescribed basis. We estimate the coefficients of the modal expansion 
using available observational data. This is usually accomplished through two distinct frameworks. 
One is deterministic and estimates the expansion coefficients by solving a least-squares (LS) 
problem. The second is probabilistic and uses a Bayesian approach to derive a distribution for 
the coefficients, resulting in the maximum-a-posteriori (MAP) estimate. Here, we seek to quantify 
and compare the accuracy of these two approaches. To this end, we derive a computable expression 
for the difference in Bayes risk between the deterministic LS and the Bayesian MAP estimates. 
We prove that this difference is always nonnegative, indicating that the MAP estimate is always 
more reliable than the LS estimate. We further show that this difference comprises two nonnegative 
components representing measurement noise and prior uncertainty, and identify regimes where 
one component dominates the other in magnitude. We also derive a novel prior distribution from 
the sample covariance matrix of the training data, and examine the greedy Bayesian and 
column-pivoted QR (CPQR) sensor placement algorithms with this prior as an input.  Using 
numerical examples, we show that the greedy Bayesian algorithm returns nearly optimal sensor 
locations. We show that, under certain conditions, the greedy Bayesian sensor locations are 
identical or nearly identical to those of CPQR when applied to a regularized modal basis. 
\end{abstract}

\begin{keyword}
state estimation \sep Bayesian inverse problems \sep least squares \sep empirical interpolation \sep sensor placement
\end{keyword}

\end{frontmatter}

\section{Introduction}

Estimating an unknown state of a system from sparse and noisy measurements
arises in a range of applications such as fluid dynamics~\cite{nishida2022,barklage2024},
weather prediction~\cite{Ghil1981,Farrell_2006}, and control theory~\cite{lin2020,barfoot2024}.
In this article, we consider linear methods which estimate the state within a modal basis. 
The observational data is used to estimate the coefficients of the expansion in this basis.

There are two common approaches to address state estimation problems:
deterministic and Bayesian.  A deterministic approach
uses a minimum-norm least squares (MNLS) estimate of the expansion coefficients, leading to the discrete empirical interpolation method (DEIM)~\cite{chatur2010,manohar2018}.  
On the other hand, the Bayesian approach
assumes that the coefficients are random variables and imposes a prior
distribution on them, which is often assumed to be Gaussian. Given the observational data, an application of
Bayes rule returns the posterior distribution, which is also Gaussian~\cite{kay1993}. It is common practice to use the
maximum-a-posteriori (MAP) estimate, which is the mode of the posterior
distribution, as a point estimator for the expansion coefficients.

This article has two main objectives.  First, we seek to bridge the
gap between the deterministic and Bayesian approaches for state estimation by
comparing the performance of the MNLS and MAP estimates theoretically and
numerically.  
We compare these two point estimates by considering their Bayes risk, which is the expected squared error. 
It is well-known that the mean of the posterior distribution minimizes the Bayes risk (see Chapter 10 of Kay \cite{kay1993}).  
When the prior distribution and measurement noise are Gaussian, the posterior distribution is Gaussian with a mean equal to its mode.
Thus, the MAP estimate is optimal in a Bayesian framework.
Here we seek to understand in what sense the MNLS estimate is suboptimal.   
We also consider the 
problem of sensor placement. 
Specifically, we consider the problem of finding an optimal subset of candidate 
sensor locations at which to measure
the state. This is done within both Bayesian and deterministic paradigms. We provide theoretical and numerical insight into how these two approaches compare.

\paragraph{Outline and contributions}

This work provides both theoretical and numerical results regarding the deterministic and Bayesian
approaches to state estimation. The specific contributions are:

\begin{enumerate}
\item Theory (discussed in \Cref{sec:theory}):
\begin{enumerate}
\item We derive an expression for the difference between the 
Bayes risk of the MNLS and MAP
estimates, which we prove to be non-negative, and identify regimes where this difference 
is large and small.  We identify two distinct components of this difference, one of 
which relates to the underlying system and the other to measurement noise.  
Moreover, we derive informative upper bounds on each component.

\item We propose a novel prior covariance matrix obtained using a change of basis
of the sample covariance matrix of the training data. 

\item We derive bounds on the relative error of the MAP estimate, which
complements the error analysis of the DEIM estimate.
\end{enumerate}
\item Sensor placement (discussed in \Cref{sec:alg}): We investigate two algorithms for sensor placement in the Bayesian case. The first is a greedy approach and the second is based on column-pivoted QR. 
These methods have been proposed and studied in the context of inverse problems previously.
Here, we adapt them to state estimation, 
and discuss the connections to deterministic approaches for sensor placement.

\item Numerical Results (discussed in \Cref{sec:num_res}):
\begin{enumerate}
\item We compare the components of the difference in Bayes risk between the DEIM and MAP estimates, 
in addition to their upper bounds.  
We provide insight into the behavior of this difference in different mode/sensor regimes.
Our numerical results confirm that both theoretical upper bounds are nearly optimal.
\item We show that as the number of modes exceeds the number of sensors, the relative
error of a greedy Bayesian algorithm remains flat or decreases, while that of a greedy deterministic
algorithm remains flat or increases.
\item We show that two sensor placement techniques return sensor locations that 
are nearly identical and nearly optimal. Furthermore, they result in reconstructions of nearly identical quality.
\end{enumerate}
\end{enumerate}

\paragraph{Related Work}
We first review the literature on deterministic methods for state estimation.
Chaturantabut and Sorensen \cite{chatur2010} first introduced the DEIM formula and DEIM
algorithm in the context of reduced-order modeling.  
Manohar et al. \cite{manohar2018} subsequently repurposed the DEIM formula for 
state estimation and sensor placement.
One drawback of the DEIM sensor placement approach is that it produces a different set of 
sensors for the same modal basis given a reordering of the basis vectors.
Such behavior is undesirable, as the ordering of basis vectors is immaterial to the 
information encoded in the basis. 
Drma\v{c} and Gugercin \cite{drmac2016} proposed an approach based on the CPQR algorithm, 
which is invariant to any reordering or rotation of the basis vectors and has improved performance bounds. In a similar approach, Drma\v{c} and Saibaba~\cite{drmac2018} use a strong rank-revealing QR factorization to dramatically improve the bounds from CPQR and removing the exponential dependence on the number of sensors. Variants of the DEIM formula 
include Sparse DEIM (S-DEIM) by Farazmand \cite{farazmand2024}, applicable to 
dynamical systems, and tensor-DEIM by Farazmand and Saibaba \cite{farazmand2023}, 
applicable to multidimensional data. Another drawback of the DEIM algorithm and formula is that neither explicitly accounts for 
measurement noise. 
Peherstorfer et al. \cite{peherstorfer2020} observe the amplification of noise-related error when the number of 
sensors equals the number of modes, which we likewise observe.  
To address this, they recommend using more sensors than modes for state estimation.

We next consider the Bayesian approaches for state estimation. 
Hirsh et al. \cite{hirsh2022} demonstrate the efficacy of sparsity-promoting priors for model uncertainty 
quantification, including the `spike and slab' and `regularized horseshoe' distributions.
Sensor placement in the Bayesian framework typically  involves optimizing some functional (e.g., trace or determinant)
of posterior covariance.
One popular solution is the greedy algorithm rigorously analyzed by 
Nemhauser and Wolsey \cite{nemhauser1978}, which provides performance guarantees for certain
types of objective functions (which we discuss in \Cref{sec:alg}).
Nishida et al. \cite{nishida2022}, Krause et al. \cite{krause2008}, and Shamaiah et al. \cite{shamaiah2010} 
have since used this algorithm to find 
sensor locations that approximately optimize a Bayesian objective.  
In our numerical results, we use this algorithm.
Alternatively, Attia et al. \cite{attia2022} treat available locations as Bernoulli random variables with 
some probability of having a sensor, allowing for stochastic optimization.  Eswar et al. \cite{eswar2024} 
frame Bayesian sensor placement as a column subset selection problem (CSSP), and propose 
randomized algorithms with high-probability lower bounds on their Bayesian objective.

A work closely related to ours is that of Klishin et al. \cite{klishin2023}, 
who decompose their Bayesian objective into terms describing sensor variance and sensor interactions.  
As with our work, one of their sensor placement algorithms
uses CPQR on the standard basis of proper orthogonal decomposition (POD) modes, 
and on the basis of POD modes regularized by a prior covariance matrix.  They use a MAP reconstruction,
and compare reconstruction error with respect to the number of sensors.
While we both assume Gaussian priors, a fundamental difference between our implementations 
is the choice of the prior covariance matrix, which we discuss in \Cref{subsec:phi_gamma_pr}.

Several works frame estimates as both deterministic and Bayesian.
Tibshirani \cite{tibshirani1996} states that, in the case of a double exponential 
(i.e. Laplace) prior distribution, the MAP estimate is the solution of a least-squares problem 
with an $L^1$-norm penalty (see Section 5 of Tibshirani \cite{tibshirani1996}).  
Gribonval \cite{gribonval2011} shows that, in the case of iid Gaussian noise,
the MAP estimate is the solution of some penalized least-squares problem for any 
prior.  To the best of our knowledge,
our work is the first to rigorously compare the Bayes risks of the MAP and DEIM 
estimates of a linear inverse problem.

In a broader context,
the present paper bridges the gap between the probabilistic Bayesian methods and the 
deterministic least-squares approach to linear inverse problems.

\section{Preliminaries}\label{sec:prelim}
\subsection{Discrete Empirical Interpolation Method (DEIM)}
We begin by describing the discrete empirical interpolation method (DEIM), first introduced by 
Chaturantabut and Sorensen~\cite{chatur2010} for reduced-order modelling and later repurposed by 
Manohar et al.~\cite{manohar2018} for state estimation.
Consider a training data matrix $X_{\text{train}} \in \R^{N \times p}$. Each column of $X_{\text{train}}$ corresponds to a snapshot (i.e.,
sample) of the system.  From $X_{\text{train}}$, we extract a basis matrix $\Phi \in \R^{N \times n}$ whose columns 
approximately span the training data. A common choice for the basis matrix is the POD modes, which coincide with the left singular vectors of the data matrix
 $X_{\text{train}}$.  Sensors that take measurements can be placed 
at any of the $N$ dimensions of the data.  The selection matrix $S \in \R^{N \times k}$ has a subset of the columns of the identity matrix as its columns and represents the locations of these sensors.  
More precisely, a sensor at the $i^{\th_}$ dimension of the data corresponds to a column of $S$ which coincides with the $i^{\th_}$ column of the $N \times N$ identity.  
Neglecting observational noise, the true state of the system $u\in\R^N$ and the vector of observations $y \in \R^k$ from the sensors are related through
\begin{equation}
    y = S^\top u. \label{eq:y_obs}
\end{equation}
Although the true state $u$ is unknown, the observations $y$ are available through sensor measurements. 

Within the orthogonal basis $\Phi$, the best approximation of $u$ is given by
\begin{equation}
    u \approx \Phi m_* \com
\end{equation}
where $m_*:=\Phi^\top u$ is the coordinate vector of the orthogonal projection of the true state $u$ onto the column space of the basis $\Phi$.  
Note that, since the true state $u$ is unknown, we cannot directly compute $m_*$.
Instead, we approximate $m_*$ by solving the least-squares problem,
\begin{equation}\label{eq:DEIM_ls_prob}
    \min_{m \in \R^n} \norm{y-S^\top \Phi m}_2,
\end{equation}
which minimizes the discrepancy between the true observations $y$ and the observations obtained from the approximation $u\approx \Phi m$.
Any solution of \Cref{eq:DEIM_ls_prob} takes the form~\cite{farazmand2024},
\begin{equation}
    \tilde{m}_{\LS}(z) = (S^\top \Phi)^+ S^\top u + z,
\end{equation}
where $z \in \R^n$ is an arbitrary vector in the null space of $S^\top \Phi$.  We use the superscript +
to denote the Moore-Penrose pseudoinverse.
The general DEIM estimate of the true state $u$ given some $z \in \text{Null}(S^\top \Phi)$ is then
\begin{align}
    \begin{split}\label{eq:deim}
        \tilde{u}_{\DEIM}(z) :=&\, \Phi \tilde{m}_{\LS}(z) \\
        =&\, \Phi (S^\top \Phi)^+ S^\top u + \Phi z \period
    \end{split}
\end{align}
The MNLS estimate of the coefficient vector $m$ (denoted herein as $m_{\LS}$) corresponds to $z=0$,
meaning $m_{\LS}:=\tilde{m}_{\LS}(0)$.  This is a particularly common estimate \cite{manohar2018, clark2021}.
We denote the MNLS estimate of the full state $u$ by $u_{\DEIM}:=\tilde{u}_{\DEIM}(0)$.
Chaturantabut and Sorensen \cite{chatur2010} used \Cref{eq:deim} only in the case of an 
invertible $S^\top \Phi$, while Manohar et al.~\cite{manohar2018} and Clark et al.~\cite{clark2021} 
considered the more general 
setting where $S^\top \Phi$ is not necessarily invertible or square.
We discuss a number of algorithms for estimating an optimal selection matrix $S$ in \Cref{sec:alg}.
\Cref{tab:notation} summarizes some notation that we use throughout this paper.

\begin{table}[H]
    \centering
    \begin{tabular}{|>{\arraybackslash}p{7cm}|>{\arraybackslash}p{6cm}|}
        \hline
        \textbf{Symbol} & \textbf{Description} \\
		\hline
        $N$ & High-fidelity resolution \\
        \hline
        $n$ & Number of modes \\
        \hline
        $k$ & Number of sensors \\
        \hline
        $S \in \R^{N \times k}$ & Selection matrix \\
        \hline
        $\Phi \in \R^{N \times n}$ & Modal basis \\
        \hline
        $A:=S^\top \Phi \in \R^{k \times n}$ & Map from the modal coordinates to observation space \\
        \hline
        $m \in \R^{n}$ & Coordinates of $u$ in the $\Phi$ basis \\
        \hline
        $m_{\LS}:=A^+y \in \R^{n}$ & Minimum-norm least squares (MNLS) estimate of $m$ \\
        \hline
        $m_{\MAP}:=\Gamma_{\post}A^\top\Gamma_{\noise}^{-1}y \in \R^{n}$ & Maximum-a-posteriori (MAP) 
        estimate of $m$ \\
        \hline
        $\Risk(\hat m)=\E_m[\E_{y|m}[\norm{\hat m(y)-m}_2^2]]$ 
        & Bayes risk of the estimate $\hat m$ of $m$ given observations $y$ 
        \\
        \hline
        $u \approx \Phi m \in \R^N$ & Unknown full state \\
		\hline
        $\eta \in \R^k$ & Measurement noise \\
        \hline
        $y =S^\top u + \eta \in \R^k$ & Observations \\
        \hline
        $\Gamma_{\text{prior}} \in \R^{n \times n}$ & Prior covariance matrix of $m$  \\
        \hline
        $\Gamma_{\text{noise}} \in \R^{k \times k}$ & Covariance matrix of $\eta$  \\
        \hline
        $\Gamma_{\text{post}}=\left(\Gamma_{\text{prior}}^{-1}+A^\top \Gamma_{\noise}^{-1}A\right)^{-1} \in \R^{n \times n}$ & Posterior covariance matrix of $m$ \\
        \hline
        $u_{\MAP}=\Phi m_{\MAP} \in \R^N$ & MAP reconstruction of $u$ \\
        \hline
        $u_{\text{DEIM}}=\Phi m_{\LS} \in \R^N$ & DEIM reconstruction of $u$
        \\
        \hline
    \end{tabular}
    \caption{Some relevant notation.}
    \label{tab:notation}
\end{table}

\subsection{The Bayesian approach}\label{subsec:overview}
We now describe the Bayesian approach to state estimation.  
Assume the following relationship between
 $m$ and $y$:
\begin{equation}
    y=Am+\eta \com \label{eq:y}
\end{equation}
where $A$ is a linear map and $\eta$ is 
measurement noise.  In our application, the linear map is $A:=S^\top \Phi$.  Our task is to estimate
the parameter $m$ from \Cref{eq:y} given measurements $y$.

Whereas DEIM uses two sources of information, namely measurements $y$ and a basis $\Phi$, 
to estimate the inversion parameter $m$, the Bayesian approach uses four: measurements $y$,
a basis $\Phi$, a probability distribution for $m$ prior
to measurements, and the distribution of measurement noise $\eta$.  The distribution of $m$ 
prior to measurements is based on expert understanding or simulation.  
The distribution of $\eta$ is based on the quality of the 
instruments taking measurements.
\begin{figure}
	\centering
	\includegraphics[width=\textwidth]{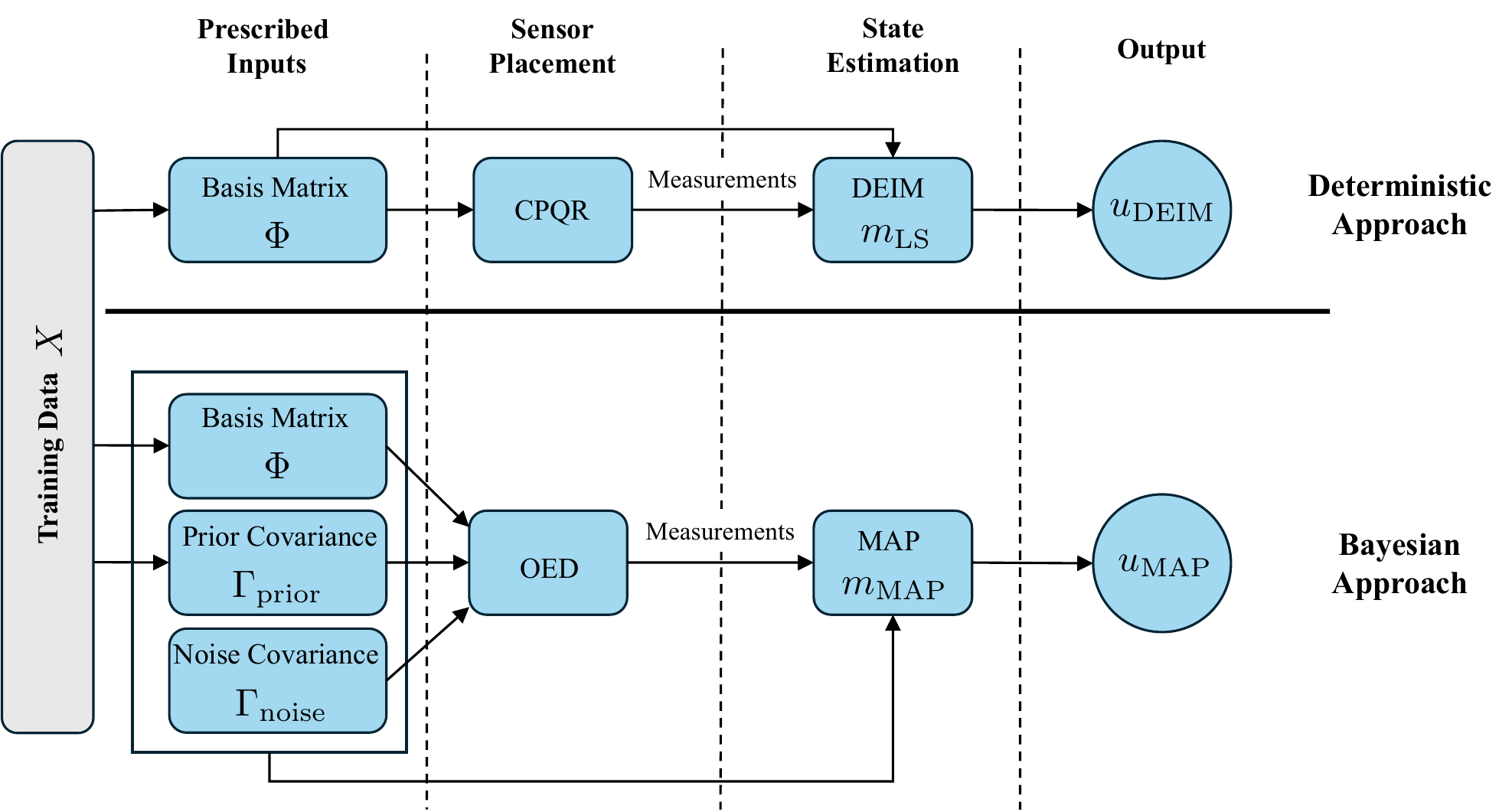}
	\caption{Summarized workflow of the deterministic and Bayesian approaches.  
		Arrows indicate the flow of information.  Prescribed inputs are the modal basis $\Phi$, 
		the prior covariance matrix $\Gamma_{\prior}$, and the noise covariance matrix $\Gamma_{\noise}$.  
		From these inputs, CPQR and the OED criterion generate the sensor locations.  
		DEIM and MAP then use measurements at these locations to 
		compute $m_{\LS}$ (the MNLS estimate of $m$) and $m_{\MAP}$ (the MAP estimate
		of $m$).  These approximations of the inversion parameter $m$ then generate $u_{\MAP}:=\Phi m_{\MAP}$ and 
		$u_{\text{DEIM}}:=\Phi m_{\LS}$, the MAP and DEIM approximations of the full-state $u$.} 
	\label{fig:workflows}
\end{figure}

Prior to taking measurements, we assume that the distribution of 
$m$ is Gaussian
\begin{equation}
    m \sim \normal(0\,, \Gamma_{\prior}), \label{eq:m}
\end{equation}
where $\Gamma_{\prior}\in \R^{n\times n}$ denotes the prior covariance matrix.
We likewise assume measurement noise is distributed as
\begin{equation}
    \eta \sim \normal(0\,, \Gamma_{\noise}), \label{eq:eta}
\end{equation}
where $\Gamma_{\noise}\in \R^{k\times k}$ is the noise covariance matrix.  
Through a direct application of Bayes' formula (see Section 3.4 of Kaipio and Somersalo~\cite{kaipio2005}), the posterior
distribution of the coordinates $m$ given measurements $y$ is 
\begin{equation}
    m|y \sim \normal(m_{\post}, \Gamma_{\post}), \label{eq:m|y}
\end{equation}
where the posterior covariance matrix is 
\begin{equation}
    \Gamma_{\post}=(\Gamma_{\prior}^{-1}+A^\top \Gamma_{\noise}^{-1}A)^{-1}, \label{eq:Gamma_post}
\end{equation}
and the posterior mean is
\begin{equation}
    m_{\post}=\Gamma_{\post}A^\top \Gamma_{\noise}^{-1}y.  \label{eq:m_post}
\end{equation}
We refer to the distribution in \Cref{eq:m} as the \emph{prior}, and that in 
\Cref{eq:m|y} as the \emph{posterior}.
Since the posterior is Gaussian, $m_{\post}$ coincides with the \emph{maximum-a-posteriori}
(MAP) estimate of $m$, which we denote by $m_{\MAP}$.  After 
taking measurements, we approximate the full state of the system as 
$u_{\MAP}:=\Phi m_{\MAP}$.

Within this framework, Bayes risk under an $L^2$ loss (or simply \emph{Bayes risk}) 
of some arbitrary estimator $\hat m$ is  
\begin{align}\label{eq:bayes_risk}
    \begin{split}
        \Risk(\hat m):&=\E_{m,y}[\norm{\hat m(y)-m}_2^2] \\
        &=\int \int \norm{\hat m(y)-m}_2^2 f(m,y) \, \id y \, \id m,
    \end{split}
\end{align}
where $f(m,y)$ is the joint probability density function of the inversion 
parameter $m$ and measurements $y$. Here, $\hat m(y)$ denotes an estimate of the inversion parameter given the observational data $y$. 
The MAP estimate~\eqref{eq:m_post} is a special choice of $\hat m$. The Bayes risk $\Risk(\hat m)$ measures the expected error in estimating the inversion parameter $m$.
We note that $\Risk(0)\geq \Risk(m_{\MAP})$, which follows as corollary of Theorem 12.1 of Kay \cite{kay1993}.
Furthermore, the Bayes risk of the MAP point coincides with the trace of the posterior 
covariance matrix (see Section 2.5 of Chaloner and Verdinelli~\cite{chaloner1995}),
\begin{equation}\label{eq:MAP_Bayes_risk}
    \Risk(m_{\MAP})=\tr(\Gamma_{\post}).
\end{equation} 
In \emph{A-optimal design}, one selects the sensor locations to minimize $\tr(\Gamma_{\post})$. Eq.~\eqref{eq:MAP_Bayes_risk} shows that A-optimal design is equivalent to minimizing the Bayes risk $\Risk(m_{\MAP})$.
On the other hand, \emph{D-optimal design} refers to sensor selection that minimizes $\log\left[\det(\Gamma_{\post})\right]$, which is equivalent to
maximizing the expected information gain.  
For an overview of these results concerning the A- and D-optimal designs, see Section 2.2 of Chaloner and Verdinelli~\cite{chaloner1995}.
The A- and D-optimal designs are two instances of optimal experimental design (OED) criteria.

To estimate the optimal sensor locations used by the MAP estimate, we use the greedy OED algorithms discussed in \Cref{sec:alg}.
\Cref{fig:workflows} summarizes the workflows of the deterministic and 
Bayesian approaches to sensor placement and state estimation.

\section{Theoretical Results}\label{sec:theory}
In this section, we present our main theoretical results involving a principled choice of the prior $\Gamma_{\prior}$ and the comparative quality of the MNLS and MAP estimates.

\subsection{Modal basis and prior covariance matrix}\label{subsec:phi_gamma_pr}
We compute the modal basis $\Phi$ and prior covariance matrix $\Gamma_{\prior}$ directly from $p$ 
snapshots of the full $N$-dimensional system. We denote these snapshots by $u_1$, $u_2$, ..., $u_p\in \R^N$ and define the data matrix,
\begin{equation}
    X:=\sqbr{(u_1 - \bar{u}) \, | \, (u_2 - \bar{u}) \, | \, ... \, | \, (u_p - \bar{u})}\in \R^{N \times p},
\end{equation}
where $\bar{u}=\frac{1}{p}\sum_{i=1}^p u_i$.
Following common practice~\cite{chatur2010,manohar2018,farazmand2024}, we derive a set of basis functions from the snapshot data using proper orthogonal decomposition.
Specifically, consider the economical singular value decomposition (SVD) of the data matrix, $X=U_{r} \Sigma_{r} V_{r}^\top $,
where $r$ is the rank of $X$.  
We take the basis matrix $\Phi=U_n$ where $U_n$ comprises the first $n$ columns of $U_{r}$ and $n\leq r$.

Given this basis, we turn our attention to a principled choice for the prior covariance matrix $\Gamma_{\prior}$.
First, consider the unbiased sample covariance matrix $\Gamma_X:=\frac{1}{p-1}XX^\top $ of the snapshots.  
Using economical SVD of the data matrix, we have
\begin{equation}
    \Gamma_X:=\frac{1}{p-1}U_{r} \Sigma_{r}^2 U_{r}^\top .
\end{equation}
On the other hand, we note that $\Gamma_{\prior}$ 
is the covariance of the coordinates of the snapshots with respect to the basis $\Phi=U_n$.  
Through a change of basis, we obtain the prior covariance matrix,
\begin{equation}\label{eq:prnciple_prior}
        \Gamma_{\prior}=U_{n}^\top   \Gamma_X U_{n} =\frac{1}{p - 1}\Sigma_n^2,
\end{equation}
where $\Sigma_n$ is a diagonal matrix containing the largest $n$ singular values of $X$.
This prior is derived systematically from the sample covariance of the snapshot data and is distinct from that of Klishin et al.~\cite{klishin2023}, who propose using either $\Sigma_n$ or a scalar multiple of the identity matrix.

\subsection{Performance guarantees of the MAP estimate}\label{subsec:map_perforance}
To assess the accuracy of a state estimate, one usually examines its relative error. The relative error is not computable for the test data where the true state $u$ is unknown.
However, a computable a priori upper bound on this error can be derived. In \Cref{subsubsec:map_rel_err}, we briefly discuss such upper bounds and argue that they are \emph{not} informative for comparing the accuracy of the deterministic and Bayesian methods. Subsequently, in~\Cref{subsubsec:map_risk_prem}, we use the Bayes risk to derive a computable and informative identity which quantifies the reliability of the deterministic method compared to the Bayesian method.
Our choice of modal basis $\Phi$ and prior covariance matrix $\Gamma_{\prior}$ has no impact on the theory presented in 
\Cref{subsubsec:map_risk_prem}, but does impact that of \Cref{subsubsec:map_rel_err}.

\subsubsection{Upper bounds on the relative error}\label{subsubsec:map_rel_err}
\Cref{prop:map_ub} presents an upper bound on the relative error of the MAP 
estimate $u_{\MAP}$ of the full state $u$.
\begin{proposition}\label{prop:map_ub}
	Let $D := \Phi \Gamma_{\post} (S^\top  \Phi)^\top  \Gamma_{\noise}^{-1} S^\top $ 
	and $\eta$ denote measurement noise.  Then
	\begin{equation}\label{eq:map_ub}
	\frac{\norm{u_{\MAP}-u}_2}{\norm{u}_2} \leq
	\norm{D-I}_2 + \norm{\Gamma_{\post} 
		(S^\top  \Phi)^\top  \Gamma_{\noise}^{-1}}_2\, 
	\frac{\norm{\eta}_2}{\norm{u}_2}.
	\end{equation}
\end{proposition}
\begin{proof}
	See \Cref{appendix:1}.
\end{proof}

For comparison, we derive an analogous upper bound on the relative 
error of the DEIM estimate $u_{\DEIM}$ of the full state $u$.
\begin{proposition}\label{prop:deim_ub}
	For the DEIM relative error, we have that
	\begin{equation}\label{eq:deim_ub}
	\frac{\norm{u_{\DEIM}-u}_2}{\norm{u}_2} \leq 
	\norm{(S^\top  \Phi)^+}_2 \,
	\paren{1+\frac{\norm{\eta}_2}{\norm{u}_2}}.
	\end{equation}
\end{proposition}
\begin{proof}
	See \Cref{appendix:1}.
\end{proof}
Although the upper bounds in~\eqref{eq:map_ub} and \eqref{eq:deim_ub} are computable, their direct comparison is not straightforward. For instance, it is not immediately clear  whether one upper bound is always smaller. 
More importantly, these upper bounds do not use the information contained in the 
prior and noise covariance matrices, and thus may be pessimistic, meaning the true relative errors can be much smaller. This motivates the use of an alternative measure for comparing the reliability of the LS and MAP estimates. In~\Cref{subsubsec:map_risk_prem}, we show that the risk premium (defined below) yields a computable and informative measure of reliability.

\subsubsection{The risk premium}\label{subsubsec:map_risk_prem}
Before we proceed, note that the theory presented in this subsection (\Cref{subsubsec:map_risk_prem}) holds
for any matrix $A\in \R^{k\times n}$, and is not restricted to $A=S^\top \Phi$. 
In this subsection, we thus use the term `MNLS estimator,' which implies a general $A\in \R^{k\times n}$, 
instead of `DEIM estimator,' which implies $A=S^\top \Phi$.

The Bayesian Gauss--Markov Theorem implies that, over all estimators of the inversion
parameter $m$ from \Cref{eq:y} that are linear in the observations $y$, the
estimator with minimum Bayes risk takes the form of~\Cref{eq:m_post}; see
Theorem 12.1 of Kay \cite{kay1993}.  We recall that, in our setting, this estimator
coincides with the MAP estimate.  

Although the MAP estimate $m_{\MAP}$ minimizes the Bayes risk, it is not immediately clear how it compares with
its deterministic counterpart $m_{\LS}$ derived from the least-squares problem~\eqref{eq:DEIM_ls_prob}.
Although $m_{\LS}$ is non-Bayesian, one can nonetheless compute its Bayes risk using~\Cref{eq:bayes_risk}. A straightforward calculation shows that
\begin{equation}\label{eq:ls_bayes_risk}
\Risk(m_{\LS})=\tr \sqbr{(I-A^+A)\Gamma_{\prior}}+
\tr\sqbr{A^+\Gamma_{\noise}(A^+)^\top }.
\end{equation}
To arrive at a fair comparison, here we used the same prior covariance $\Gamma_{\prior}$ as for the MAP estimate.
In \Cref{eq:ls_bayes_risk}, we see that the Bayes risk of the MNLS estimator consists of two distinct 
components: that which stems from prior uncertainty of the inversion parameter and that
which stems from measurement noise.

In~\Cref{thm:bayes_risk_ineq}, we derive a computable and informative expression for the difference $\Risk(m_{\LS})-\Risk(m_{\MAP})$, which we call 
the risk premium. Before stating our results, we define a number of quantities which appear in this theorem.
First, we define
\begin{align}
    \Delta_{\prior}&:= (I-A^+A)\paren{\Gamma_{\prior}-\Gamma_{\post}} \in \R^{n \times n}, \label{eq:Delta_prior}
\end{align}
which is the orthogonal projection of
$\Gamma_{\prior}-\Gamma_{\post}$ onto the null space of $A$.  
We also define
\begin{equation}\label{eq:Delta_noise}
    \Delta_{\noise} := A^+ \Gamma_{\noise} (A^+)^\top - A^+A \Gamma_{\post} \in \R^{n \times n},
\end{equation}
whose range is a subset of the range of $A^\top$.  By the fundamental theorem of linear algebra,
we have
\begin{equation}
    \Delta_{\prior}^\top \Delta_{\noise}=0_{n \times n}.
\end{equation}

We denote the traces of these matrices by
\begin{align}
	\delta_{\prior}&:=\tr \paren{\Delta_{\prior}},\label{eq:delta_prior} \\
    \delta_{\noise}&:=\tr \paren{\Delta_{\noise}}.\label{eq:delta_noise}
\end{align}

The quantities $\delta_{\prior}$ and $\delta_{\noise}$ allow for a precise comparison of
the statistical properties of the MAP and MNLS estimates.  Each estimate has an associated 
Bayes risk, which represents the expected squared error of that estimate.  
Intuitively, the MAP estimate, which accounts for the prior and measurement noise distributions, 
should have a lower Bayes risk than the MNLS estimate.  To quantify the difference in Bayes risks between these 
estimates, we introduce the following definition.

\begin{definition}[Risk Premium]\label{def:risk_prem}
    In the context of the linear Bayesian inverse problem~\Cref{eq:y}, we define the \emph{risk premium} between the
    $m_{\LS}$ and $m_{\MAP}$ estimates as $\Risk(m_{\LS})-\Risk(m_{\MAP})$.
\end{definition}

The following theorem quantifies the risk premium between the $m_{\LS}$ and $m_{\MAP}$ in terms of $\delta_{\prior}$ and $\delta_{\noise}$
and subsequently shows that the risk premium is always non-negative.

\begin{theorem}\label{thm:bayes_risk_ineq}
    Let $y=Am+\eta$, where $m \sim N(0,\Gamma_{\prior})$, 
    $\eta \sim N(0, \Gamma_{\noise})$, and $A \in \R^{k \times n}$.  
    Let $m_{\LS}:=A^+y$ denote the MNLS estimate 
    of the inversion parameter $m$ given observations $y$ 
    and $m_{\MAP}:=\Gamma_{\post}A^\top\Gamma_{\noise}^{-1}y$ 
    denote the MAP estimate.  Then, given the definitions of 
    $\delta_{\prior}$ in \Cref{eq:delta_prior} and 
    $\delta_{\noise}$ in \Cref{eq:delta_noise}, 
    we have that
    \begin{align}
         \delta_{\prior} &\geq 0, \\
         \delta_{\noise} &\geq 0, \\
         \Risk(m_{\LS}) - \Risk(m_{\MAP}) &= \delta_{\prior} + \delta_{\noise}.
    \end{align}
\end{theorem}
\begin{proof}
    See \Cref{appendix:2}.
\end{proof}

As \Cref{eq:y} shows, variance of the inversion parameter $m$ in the null space of $A$ does not result
in variance in the measurements $y$.  Variance of 
$m$ in the range of $A^\top$ does result in variance in $y$.
We thus refer to variance of $m$ in the null space of $A$ as `invisible' in $y$,
and variance of $m$ in the range of $A^\top$ as `visible' in $y$.
The statement $\delta_{\prior} \geq 0$ implies that the magnitude of the posterior covariance
that is invisible in $y$ (i.e. $\tr \sqbr{(I-A^+ A)\,\Gamma_{\post}}$) is 
never greater than the magnitude of the prior covariance that is invisible in
$y$ (i.e. $\tr \sqbr{(I-A^+ A)\, \Gamma_{\prior}}$).  Likewise, 
$\delta_{\noise} \geq 0$ implies that
the magnitude of the posterior covariance that is visible in $y$ 
(i.e. $A^+ A\, \Gamma_{\post}$) is never greater than the magnitude of the 
covariance of the least-squares estimate with respect to noise 
(i.e. $\tr \sqbr{A^+ \Gamma_{\noise} (A^+)^\top}$).
We interpret the risk premium from \Cref{def:risk_prem} as the reduction in Bayes risk
associated with using a regularized (as opposed to unregularized) estimate.
Intuitively, $\delta_{\prior}$ represents the 
regularization-induced reduction in Bayes risk stemming from prior uncertainty, 
while $\delta_{\noise}$ represents the regularization-induced reduction 
in Bayes risk stemming from noise uncertainty.

In the following result, we derive an informative upper bound on 
the risk premium by bounding its constituent parts. Here, and henceforth, the eigenvalues of a symmetric matrix $M \in \R^{n\times n}$ are 
denoted by $\lambda_i(M)$ for $1 \le i \le n$, and are listed in decreasing order.

\begin{theorem}\label{thm:risk_prem_ub}
The risk premium is bounded as
    \begin{equation}\label{eq:risk_prem_ub}
            0 \leq \Risk(m_{\LS}) - \Risk(m_{\MAP}) 
            \leq \tr \sqbr{A^+ \Gamma_{\noise} (A^+)^\top} +
            \sum_{i=1}^{\textnormal{Nullity}\paren{A}} 
            \lambda_i(\Gamma_{\prior}-\Gamma_{\post}).
    \end{equation}
\end{theorem}
\begin{proof}
    See \Cref{appendix:3}.
\end{proof}

This result quantifies the extent to which the MAP estimate $m_{\MAP}$ can be more reliable than the MNLS estimate $m_{\LS}$. Regarding the nullity of $A \in \R^{k\times n}$ that appears in \Cref{eq:risk_prem_ub}, 
note that $\text{Nullity} (A)=n-\text{Rank} (A)$ by the rank-nullity theorem, with $\text{Rank} (A) \leq \min\{k,n\}.$

We denote the components of the upper bound in~\Cref{eq:risk_prem_ub} by
\begin{align}
    \zeta_{\noise} :&= \tr \sqbr{A^+ \Gamma_{\noise} (A^+)^\top}, \\
    \zeta_{\prior} :&= \sum_{i=1}^{\textnormal{Nullity}\paren{A}} 
    \lambda_i(\Gamma_{\prior}-\Gamma_{\post}).
\end{align}
As is shown in the proof of \Cref{thm:risk_prem_ub}, $\zeta_{\noise}$
and $\zeta_{\prior}$ are upper bounds on $\delta_{\noise}$ and $\delta_{\prior}$, respectively.
Large eigenvalues of $\Gamma_{\prior} - \Gamma_{\post}$ correspond to large values of $\zeta_{\prior}$,
while eigenvectors of $\Gamma_{\prior} - \Gamma_{\post}$ being close to the null space of $A$ corresponds
to $\delta_{\prior}$ being close to its upper bound $\zeta_{\prior}$.
Increasing the nullity of $A$ by adding a column to the basis $\Phi$ does not necessarily correspond to an 
increase in $\zeta_{\prior}$, as the eigenvalues of $\Gamma_{\prior}-\Gamma_{\post}$ can decrease in magnitude.

When noise is independent and identically distributed (iid), the $\zeta_{\noise}$ term simplifies to
\begin{equation}\label{eq:zeta_noise_iid}
    \zeta_{\noise}=\sigma_{\noise}^2 \, \sum_{i=1}^{\Rank{A}}\frac{1}{\sigma_i \paren{A}^2},
\end{equation}
where $\sigma_i(A)$ is the $i^{\th_}$ largest singular value of $A$ and $\sigma_{\noise}$ is the noise standard deviation.

Finally, we note that although noise does not directly enter the calculation of the MNLS estimate $m_{\LS}$, noisy data naturally affects the quality of this estimate. The Bayes risk of MNLS~\eqref{eq:ls_bayes_risk} reveals this dependence through the second term depending on $\Gamma_{\noise}$.
From \Cref{eq:ls_bayes_risk}, we also see that the $\zeta_{\noise}$ term in the upper bound in~\Cref{eq:risk_prem_ub} 
coincides with the noise-dependent term of $\Risk(m_{\LS})$.

\section{Sensor Placement Algorithms}\label{sec:alg}
The quality of reconstructions is significantly impacted by the position of the sensors that gather the observational data.
Thus, we turn our discussion to the sensor placement problem in DEIM (Section~\ref{ssec:qdeim}) and the Bayesian methods (Sections~\ref{subsec:bayesian_sensors} and~\ref{subsec:qmap}).  
In the context of Bayesian approaches, we focus our discussion on the D-optimal design as opposed to the A-optimal design.
Our reason for doing so stems from our use of a greedy algorithm for Bayesian sensor placement.
Under the assumption of iid Gaussian noise, this algorithm is guaranteed to approximately optimize the D-optimal objective, as we discuss in \Cref{subsec:bayesian_sensors}.
Such a guarantee does not exist in the case of A-optimality (see Section 8.4 of Krause et al. \cite{krause2008}).

Throughout this section, we denote by $S \in \R^{N\times k}$ a selection matrix that contains $k$ distinct columns from the $N \times N$ identity matrix.
We denote the set of $k$ sensor locations (i.e. indices) as 
$\{ \xi_1,\,\xi_2,\,...\,,\xi_k \}\subset \{ 1,\,2,\,...,\,N \}$.
These indices uniquely determine the selection matrix $S$ as
\begin{equation}
    S:=\sqbr{e_{\xi_1}\,e_{\xi_2}\,...\,e_{\xi_k}},
\end{equation}
where $e_{\xi_j}$ is the $\xi_j^{\th_}$ column of the $N\times N$ identity matrix. Recall that the observations are related to the true state by $y = S^\top u+\eta$ where $\eta$ is the measurement noise.
We assume henceforth that the noise covariance matrix, $\Gamma_{\noise} = \sigma_{\noise}^2 I$, is a multiple of the identity.

\subsection{Q-DEIM}\label{ssec:qdeim}
In the context of the DEIM state estimate \Cref{eq:deim},
a popular approach for determining the selection matrix $S$ is the column-pivoted QR 
(CPQR) algorithm proposed by Drma\v{c} and Gugercin \cite{drmac2016}. Recall the relative error upper bound~\Cref{eq:deim_ub} which involves the matrix norm $\|(S^\top \Phi)^+\|_2$.
CPQR is a greedy algorithm that attempts to minimize this matrix norm, hence minimizing the upper bound on the relative error. Note that $\|(S^\top \Phi)^+\|_2$ coincides with the largest singular value of $(S^\top \Phi)^+$, or equivalently the smallest non-zero singular value of $S^\top \Phi$. Therefore, minimizing $\|(S^\top \Phi)^+\|_2$ is equivalent to maximizing the smallest non-zero singular value of $S^\top \Phi$.

We now discuss CPQR in detail. In the discussion below, we first assume that $k \le n$; that is, the number of selected sensors is at most the number of modes. Consider the pivoted QR factorization of $\Phi^\top$,
\begin{equation}\label{eq:pivotedqr} 
\Phi^\top \bmat{\Pi_1 & \Pi_2 }  = \bmat{Q_1 & Q_2} \bmat{R_{11} & R_{12} \\ & R_{22}},
\end{equation}
where $\Pi = [\Pi_1\ \Pi_2]$ is a permutation matrix, $\Pi_1, Q_1 \in \R^{N\times k}$, and $R_{11} \in \R^{k\times k}$. The size of the other matrices should be clear from the context. The selection operator $S$ is then set equal to $\Pi_1$. The analysis in Gu and Eisenstat \cite[Theorem 7.2]{gu1996efficient} shows that 
\[ \sigma_i(R_{11}) \geq \frac{\sigma_i(\Phi)}{\sqrt{n-i}\cdot 2^i}, \qquad 1 \le i \le k.\]
From $\Phi^\top S = Q_1 R_{11}$, and the fact that $Q_1$ has orthonormal columns, we have 
\[\|(S^\top \Phi)^+\|_2= \|(\Phi^\top S)^+\|_2  = \| R_{11}^{-1}\|_2\leq \sqrt{N-k}\cdot 2^k . \]
Note here, we have used the fact that $\Phi$ has orthonormal columns so all its singular values are $1$.  If $k > n$, then once again we have the pivoted QR factorization,
\[ \Phi^\top \bmat{\Pi_1 & \Pi_2 }  = {Q_1} \bmat{R_{11} & R_{12}},\]
where $\Pi_1 \in \R^{N\times n}$, and $Q_1, R_{11} \in \R^{n\times n}$. Now, defining the permutation matrix $\Pi= \bmat{\Pi_1 & \Pi_2}$, we take $S$ to be the first $k$ columns of $\Pi$. Then, by Gu and Eisenstat \cite[Theorem 7.2]{gu1996efficient}, 
\[ \|(S^\top \Phi)^+\|_2=\|(\Phi^\top S)^+\|_2 \le \|(\Phi^\top \Pi_1)^{-1}\|_2 \le \sqrt{N-n}\cdot 2^n.\]
The first inequality follows from interlacing property of singular values, see e.g., Golub and Van Loan \cite[Corollary 8.6.3.]{golub2012matrix}.  
The above analysis generalizes the result of Drma\v{c} and Gugercin \cite[Theorem 2.1]{drmac2016}, which was done for 
the case where number of sensors $k$ equals the number of modes $n$, to the case where $k \neq n$.

 As can be seen from these bounds, in the worst case there is an exponential dependence on the number of sensors. In practice, however, CPQR performs remarkably well, and it appears that outside of some contrived test matrices this exponential dependence is not seen. Golub and Van Loan \cite{golub2012matrix} show that the cost of CPQR is $4Nnk - 2k^2 (N + n) + 4k^3$ floating point operations (flops).

 The strongest theoretical guarantees are available if we use strong rank-revealing QR (sRRQR) \cite[Algorithm 4]{gu1996efficient} instead of CPQR. In this case, by Gu and Eisenstat \cite[Theorem 3.2]{gu1996efficient}, we can obtain the analogous bounds for $k \le n$ as
\[ \|(S^\top \Phi)^+\|_2=\|(\Phi^\top S)^+\|_2  = \| R_{11}^{-1}\|_2\leq \sqrt{1 + f^2k(N-k)}, \]
where $f \geq 1$ is a user-defined factor in sRRQR. However, the computational cost of sRRQR is higher than CPQR. Therefore, in our numerical experiments we use CPQR in conjunction with the DEIM formula~\Cref{eq:deim}. This will be referred to as \emph{Q-DEIM}.

\subsection{Greedy Bayesian sensor placement}\label{subsec:bayesian_sensors}
Let $\Theta(S)$ denote an objective function which depends on the selection matrix $S$ to be minimized. 
When using the D-optimal criterion, this objective is
\begin{equation}
    \Theta_D(S):=\log \det \sqbr{\Gamma_{\post}(S)},
\end{equation}
where the posterior covariance matrix $\Gamma_{\post}(S)$ that depends on the selection operator $S$ is 
\begin{equation}\label{eq:Gamma_post__func}
    \Gamma_{\post}(S)=\left(\Gamma_{\prior}^{-1}+\sigma_{\noise}^{-2}(S^\top \Phi)^\top 
    (S^\top \Phi)\right)^{-1}.
\end{equation}
As discussed in \Cref{subsec:overview},
minimizing this objective function corresponds to maximizing the expected information gain. 
To simplify notation, we also overload the symbol $\Theta$ to represent a function of the set of sensor locations,
\begin{equation}
    \Theta(\{j_1,j_2,...,j_k\}):=\Theta(S_{\{ j_1,j_2,...,j_k \}}),
\end{equation}
where $S_{\{ j_1,j_2,...,j_k \}}$ is the selection matrix 
corresponding to the set of sensor locations $\{j_1,j_2,...,j_k\}$. If no observations are available, we have $\Theta_D(\emptyset)=\log \deter{\Gamma_{\prior}}$, where $\emptyset$ denotes the empty set.

To minimize or approximately minimize $\Theta$ with $k$ sensors, 
we use two algorithms.  The first is a brute-force search of all 
possible sensor permutations, 
which requires $\binom{N}{k}$ evaluations of $\Theta$.  Note that  
computation of $\Theta$ has a time complexity of $\Ord{n^3}$.  
For a large dimension $N$ and number of sensors $k$, this approach 
is computationally infeasible.

When faced with this limitation, we use a greedy approach to sensor 
placement described by 
\cite{nishida2022,nemhauser1978,krause2008,shamaiah2010}
and summarized in \Cref{alg:greedy_oed}.  This algorithm iterates 
$k$ times, placing a sensor at each iteration.  
At the $i^{\text{th}}$ iteration, the algorithm 
evaluates the objective gain at each of the $N-i+1$ locations 
that do not yet have a sensor.
The objective gain of a new sensor location $\ell$,
given sensors placed at locations $j_1, j_2, ..., j_{i-1}$, is
\begin{equation}\label{eq:info_gain}
    \Theta(\{ j_1, j_2, ..., j_{i-1} \})
    -\Theta(\{ j_1, j_2, ..., j_{i-1} \} \cup \{ \ell \})\period
\end{equation}
Whichever available location maximizes \Cref{eq:info_gain} 
receives the $i^{\text{th}}$ sensor.
With $N-i+1$ evaluations of 
\Cref{eq:info_gain} at the $i^{\text{th}}$ iteration
and $k$ total iterations,
the algorithm requires $\Ord{Nk}$ total evaluations of 
\Cref{eq:info_gain}.  The greedy algorithm is thus feasible even when
the dimension $N$ and number of sensors $k$ are large enough to render 
brute-force sensor placement infeasible.

\begin{algorithm}[!ht]
    \caption{The greedy approach}
    \label{alg:greedy_oed}
    \KwIn{objective function 
    $\Theta(\{j_1,j_2,...,j_i\})$
    \newline number of available sensors $k$
    }
    \KwResult{sensor locations $\xi=\{ \xi_1, \xi_2, ..., \xi_k \}$}
    
    $\mathcal{C} \gets \{ 1, 2, ..., N \}$ \\
    $\xi \gets \{ \}$ \\
    \For{$i=1,2,...,k$}{
        $\xi_i \gets \arg 
        \underset{\ell \in \mathcal{C \backslash \xi}}{\max}\left[
        \Theta(\xi)-\Theta(\xi \cup \{ \ell \})\right]$ \\
        $\xi\gets \xi\cup\{\xi_i\}$
    }
\end{algorithm}

When using the D-optimal objective $\Theta_D$, Shamaiah et al.
\cite{shamaiah2010} show a more computationally efficient 
version of \Cref{alg:greedy_oed}.
In this case, they show that a rank-1 update of 
matrix inverses allows for \Cref{eq:info_gain} to be computed at 
a cost of $\Ord{n^2}$ (see Section 4.B. of Shamaiah et al. \cite{shamaiah2010}).  
Using their rank-1 update, the flop count of \Cref{alg:greedy_oed} is 
$\Ord{Nkn^2}$, which is still higher than the $\Ord{Nkn}$ flop count of CPQR.

\Cref{alg:greedy_oed} also has performance
guarantees depending on the properties of $\Theta$.
We define the objective gain,
\begin{equation}\label{eq:obj_gain}
    J(\xi):=\Theta(\emptyset)-\Theta(\xi),
\end{equation}
of a set of sensor locations $\xi$. In the case of the D-optimal criterion, we have
$J_D(\xi):=\Theta_D(\emptyset)-\Theta_D(\xi).$
\Cref{alg:greedy_oed} attempts to approximately maximize $J$.
If objective gain~\Cref{eq:obj_gain} is monotone  
and submodular, then the set $\xi$ produced by \Cref{alg:greedy_oed}
always satisfies
\begin{equation}\label{eq:greedy_gurantee}
    \paren{1-\frac{1}{e}}J(\xi_{\opt}) \leq J(\xi)\leq J(\xi_{\opt}),
\end{equation}
where $\xi_{\opt}$ is a set of size $k$ which maximizes $J$. 
This bound follows from  
Theorem 4.2 of Nemhauser and Wolsey \cite{nemhauser1978} and shows that the sensor locations obtained from the greedy algorithm are away from the optimal sensor locations by at most a factor of $1-e^{-1}\simeq 0.63$.
Note that the D-optimal gain $J_D$ is not generally monotone and submodular (see Remark 14 of Krause et al. \cite{krause2008}). 
However, by our assumption $\Gamma_{\noise}=\sigma_{\noise}^2 I$, $J_D$ is both monotone and submodular 
(see Lemma 1 of Shamaiah et al. \cite{shamaiah2010}).
More specifically, in the case of iid Gaussian noise, we have that
\begin{equation}\label{eq:dopt_greedy_gurantee} 
    \paren{1-\frac{1}{e}}J_D(\xi_{\Dopt}) \leq
    J_D(\tilde{\xi}_{\Dopt}) \leq J_D(\xi_{\Dopt}) \com
\end{equation}
where $\xi_{\Dopt}$ is the D-optimal set of $k$ sensors,
and $\tilde{\xi}_{\Dopt}$ is the approximately D-optimal set of 
$k$ sensors produced by \Cref{alg:greedy_oed}.

Since the assumption of iid Gaussian noise is reasonable in many 
applications and, in the case
of the D-optimal criterion, results in the performance guarantee 
\Cref{eq:dopt_greedy_gurantee},
we use \Cref{alg:greedy_oed} with the objective $\Theta_D$ to select sensors.

\subsection{Q-MAP}\label{subsec:qmap}
In addition to CPQR and \Cref{alg:greedy_oed}, we consider a variant of CPQR
that uses $\Gamma_{\prior}$ in addition to $\Phi$ to place sensors.  
To motivate this algorithm, we consider the D-optimal information gain,
\begin{equation}\label{eq:Dobj_alt}
    J_D(S) = \Theta_D(\emptyset) - \Theta_D(S)  =
    \log \sqbr{ \det \paren{I + FS(FS)^\top }},
\end{equation}
where $\Theta_D(\emptyset) = \log\sqbr{\det\paren{\Gamma_{\prior}}}$ and $F:=\sigma_{\noise}^{-1}\Gamma_{\prior}^{1/2} \Phi^\top $. 
As before, we seek a selection matrix $S$ that maximizes the information gain $J_D$.

The key insight of Eswar et al. \cite{eswar2024} is that this optimization problem is closely related to column subset selection on the matrix $F \in \R^{n\times N}$. With this observation, we can perform CPQR on $F$ rather than $\Phi^\top $. This results in the pivoted QR factorization (as in~\eqref{eq:pivotedqr})
\[  F\bmat{\Pi_1 & \Pi_2} = \bmat{Q_1 & Q_2} \bmat{R_{11} & R_{12} \\ & R_{22}}.\]
The selected matrix is then given by $S = \Pi_1$. We refer to the strategy that combines
this sensor placement algorithm with a MAP reconstruction as \emph{Q-MAP}.  

Let us also denote by $\Sigma_k \in \R^{k \times k}$ a diagonal matrix whose diagonals contain the dominant singular values of $F$. By a similar argument to Eswar et al. \cite{eswar2024}, it can be shown that 
\begin{equation}\label{eq:doptbound} \logdet \paren{I +\frac{1}{q(N,k)^2}\Sigma_k^2 }  \le J_D(S) \le J_D(S_{\rm opt}) \le \logdet \paren{I + \Sigma_k^2 } , \end{equation}
where $q(N,k) = \sqrt{N-k}\cdot 2^k$ and $S_{\rm opt}$ denotes an optimal selection matrix. If sRRQR is used instead of CPQR, the bound in~\eqref{eq:doptbound} can be improved with $q(N,k) = \sqrt{1+f^2k(N-k)}$ (see the discussion in Section~\ref{ssec:qdeim}). As can be seen from~\eqref{eq:doptbound}, the set of sensors determined by the selection matrix $S$ is nearly optimal if sRRQR is used. However, in practice we use CPQR since it has a lower computational cost and excellent performance.

\section{Numerical Results}\label{sec:num_res}
We use two test problems to evaluate the sensor placement and state estimation 
approaches under study. The first problem involves the reconstruction of harmonic functions with random phases and amplitudes.  
The second problem is the reconstruction of a turbulent fluid flow 
in two dimensions~\cite{farazmand2023}.
With these test problems, we compare the greedy D-optimal sensors to those of
alternative algorithms, and the MAP estimate to MNLS.
We further compute the risk premium components ($\delta_{\noise}, \delta_{\prior}$) and their respective 
upper bounds ($\zeta_{\noise}, \zeta_{\prior}$).

In our results, we pair deterministic sensor placement algorithms with 
deterministic state estimation, and Bayesian sensor placement 
algorithms with Bayesian state estimation.  We refer to the MAP estimate
computed on the D-optimal sensor locations (found via brute-force search) as \emph{D-MAP},
and the MAP estimate computed on the greedy D-optimal sensor locations 
(found via \Cref{alg:greedy_oed}) as \emph{greedy D-MAP}.
\Cref{tab:pairings} summarizes the different pairs of sensor placement algorithms 
and state estimation formulas that we consider.
Note that in the case of Q-DEIM, we place sensors using CPQR 
with the basis $\Phi^\top$ as an input, while in the case of Q-MAP we use 
CPQR with the regularized basis $\Gamma_{\prior}^{1/2} \Phi^\top$
as an input.  In the case of Q-MAP we disregard $\Gamma_{\noise}$ since it is a scalar 
multiple of identity, and thus has no impact on the output of the CPQR algorithm.

\begin{table}[!ht]
    \centering
    \begin{tabular}{|c|c|c|c|}
        \hline
        \backslashbox{State \\ Estimation}{Sensor \\ Placement} & CPQR & D-opt. & Greedy D-opt. \\ \hline
        DEIM & Q-DEIM & ------ & ------ \\ \hline
        MAP & Q-MAP & D-MAP & Greedy D-MAP \\ \hline
    \end{tabular}
    \caption{The different pairs of sensor placement algorithms and state estimation formulas 
    that we consider.  We do not pair DEIM with the Bayesian sensor placement algorithms.}
    \label{tab:pairings}
\end{table}

\subsection{Benchmark: Random harmonic functions}\label{subsec:fourier}
We generate a data set of harmonic functions,
\begin{equation}\label{eq:fourier_data}
    f_i(x)=\sum_{j=1}^{J} \left[ a_{ij} 
    \, \sin{(j \, x + \phi_{ij})} \right],
\end{equation}
where $\phi_{ij}$ are iid random phases drawn from the uniform distribution $U[0,2\pi]$ and amplitudes
\begin{equation}\label{eq:ampl}
    a_{ij} \sim \begin{cases}
        \normal(0, 1/j), & j=1,2,...,10 \\
        \normal(0, 1/j^3), & j=11,12,...,J.
    \end{cases}
\end{equation}
The piecewise definition of amplitude given in \Cref{eq:ampl} serves 
to introduce the first spectral gap in the data after the $20^{\text{th}}$ POD mode.  
The $i^{\text{th}}$ sample is a single realization of the random function 
$f_i: [0,2\pi] \rightarrow \R$, evaluated over the spatial grid
$x=\sqbr{0, (2\pi/N), ..., (N-1)(2\pi/N)}^\top$, where $N=40$ is the high-fidelity resolution.
The number of terms in~\Cref{eq:fourier_data} is set to $J=20$. 

We generate one thousand samples of the harmonic functions where 750 are used for training and the remaining 250 are used for testing.  
We use the training data to compute the modal basis $\Phi$, the prior covariance matrix 
$\Gamma_{\prior}$, and the Bayes risk. We use the test data to compute the reconstruction error.  
We add iid Gaussian noise to the test data with mean $0$ and 
standard deviation $\sigma_{\noise}=0.1$ so that the noise covariance matrix is given by $\Gamma_{\noise}=\sigma_{\noise}^2 I$.  
On average, noise represents approximately 14.5\% of a test data sample.  Such high 
noise ensures that the effect of noise-related regularization is significant.

\subsubsection{Reconstruction error and the risk premium}\label{subsubsec:fourier_recon}
\begin{figure}[!ht]
    \hspace*{-.5cm}
    \centering
    \includegraphics[width=\textwidth]{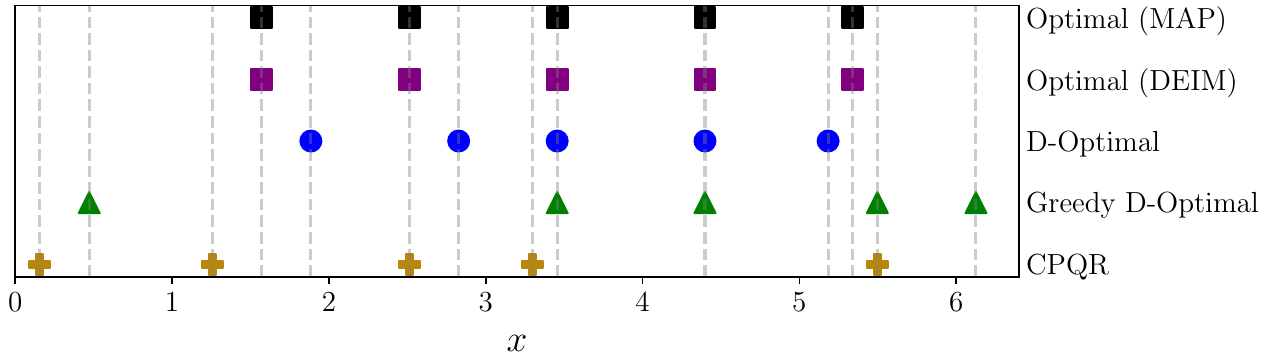}
    \caption{Sensor placement for the harmonic functions.  
    The algorithms place 5 sensors into 40 available 
    locations.  Brute force determines the optimal and D-optimal 
    sensor locations, while \Cref{alg:greedy_oed} determines the greedy D-optimal 
    sensor locations.  The optimal MAP sensor locations
	minimize the MAP reconstruction error with $k=5$ sensors and $n=20$ modes,
	while the optimal DEIM sensor locations minimize the DEIM reconstruction error
	with $k=n=5$.  We perform CPQR on the first $k$ POD modes 
    (i.e. the first $k$ columns of $\Phi$).}
    \label{fig:fourier_loc}
\end{figure}

\begin{figure}[!ht]
    \centering
    \vspace*{-2.5cm}
    \hspace*{-1cm}
    \includegraphics[width=\textwidth]{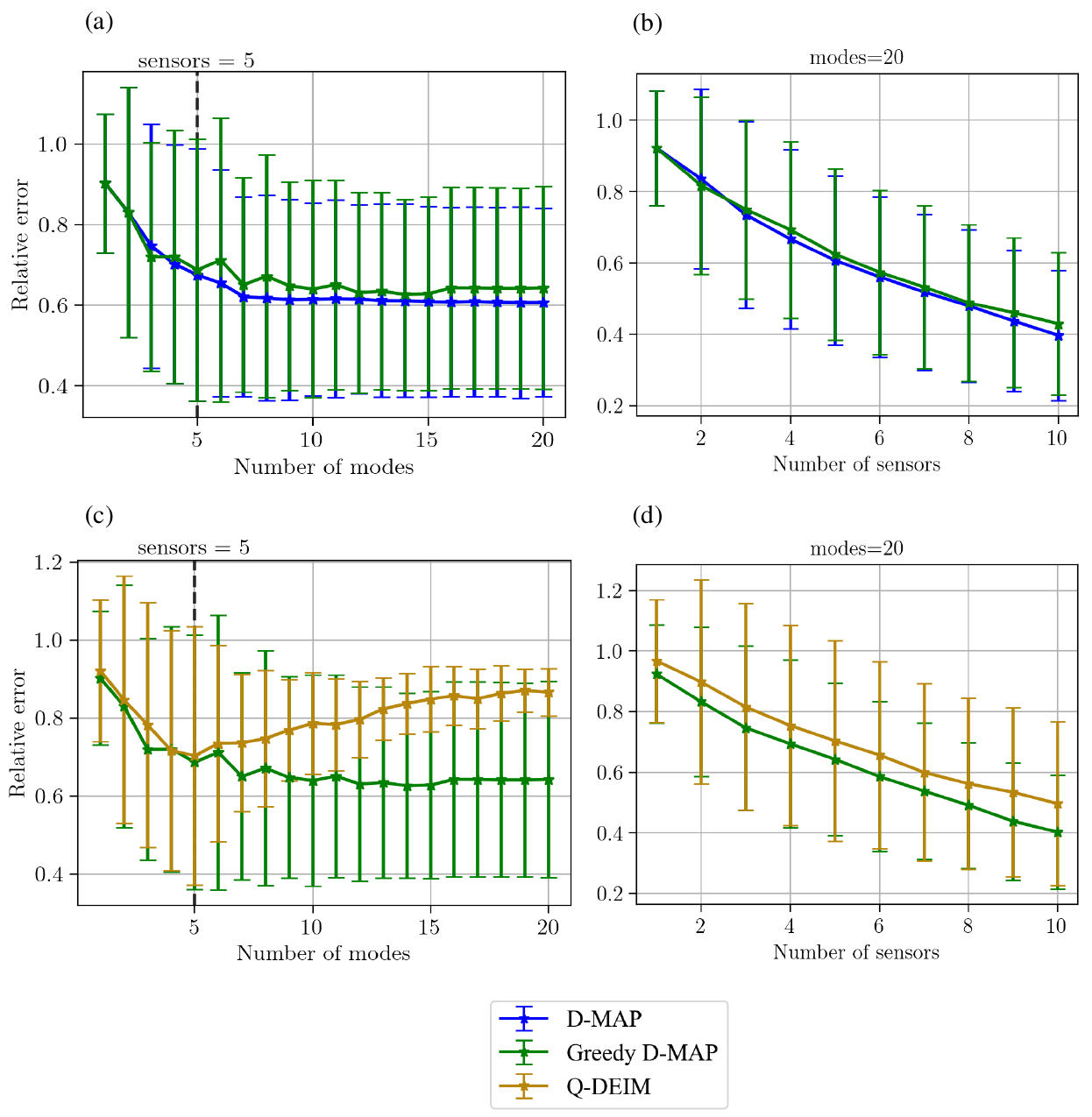}
    \caption{The relative error of D-MAP and greedy D-MAP
    on the harmonic data with respect to number of modes in (a) and sensors in (b).  
    We plot the relative 
    error of greedy D-MAP and Q-DEIM with respect to number of modes in 
    (c) and sensors in (d).
    `Relative error' refers to average relative error across all test samples, with
    error bars indicating one standard deviation.
    The `number of modes' refers to modes used for sensor placement 
    and state estimation.
    In (b), we set the high fidelity resolution $N$ to 20 for computational 
    tractability, while $N=40$ in (a), (c), and (d).  In (d), we set the number of
    modes used by Q-DEIM to the number of sensors (on the horizontal axis), and
    the number of modes used by greedy D-MAP to 20.
}
    \label{fig:Fourier_errors}
\end{figure}

\Cref{fig:fourier_loc} compares the locations of $k=5$ sensors obtained by the D-optimal, greedy D-optimal,
and CPQR algorithms to the truly optimal sensor locations. Over the space of all possible
sensor permutations, we determine two optimal permutations: 
the one that minimizes the reconstruction error of the MAP estimate on the test 
data, and the one that minimizes the reconstruction error of the DEIM estimate 
on the test data.
Given the relatively low resolution $N=40$ and only $k=5$ sensors to place,  
finding the optimal permutations is computationally feasible via a brute force search.  
In accordance with the location of the first spectral gap, we set
the number of POD modes to $n=20$ for D-MAP and greedy D-MAP.
Unlike the D-optimal and greedy D-optimal sensor locations, the CPQR
sensor locations are of highest quality when the number of modes $n$ equals the number of 
sensors $k$ (as we later see in \Cref{fig:Fourier_errors}).  
We thus set $n=k$ to generate the CPQR sensor locations.

In case of \Cref{fig:fourier_loc}, the sensor locations that minimize 
the MAP reconstruction error also minimize the DEIM reconstruction error.
We observe that the D-optimal and greedy D-optimal
algorithms select two of the optimal sensor locations, whereas the CPQR
algorithm selects one.  All the D-optimal sensors
are close to the optimal locations, whereas the greedy D-optimal and CPQR algorithms
select some locations that are further away near $x=0$ and $2\pi$.
Note that, since the harmonic functions are periodic in $x$, locations on the far left of the domain 
are close to locations on the far right.

\begin{figure}[!ht]
    \hspace*{-.9cm}
    \centering
    \includegraphics[width=\textwidth]{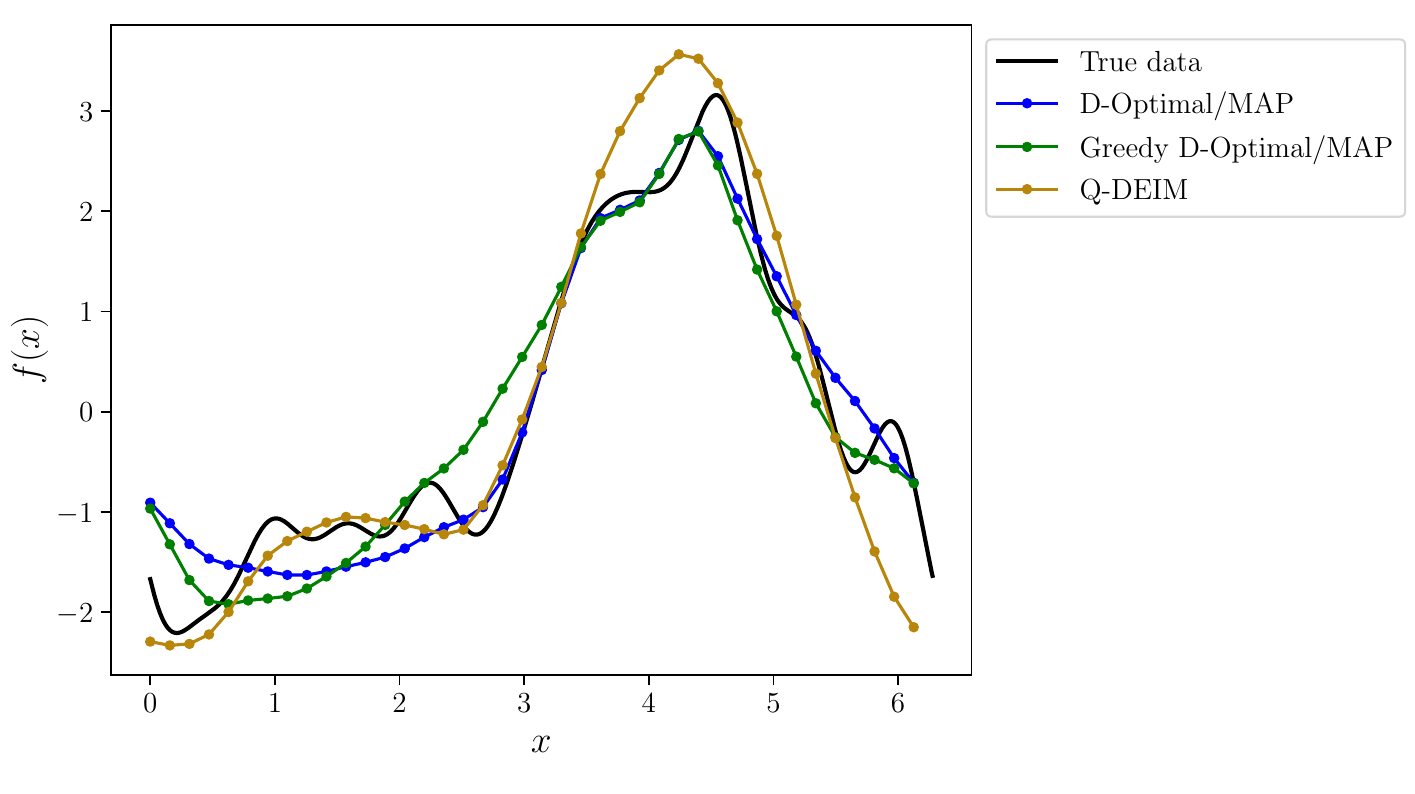}
    \caption{The reconstructions of a single sample from the harmonic functions.  
    The `true data' to be reconstructed is the black curve.
    We use 5 modes for the DEIM reconstruction (such that the number of 
    modes equals the number of sensors), and 20
    modes for the MAP reconstruction.  
    On this sample, D-MAP, greedy D-MAP, and 
    Q-DEIM have relative errors of 24.68\%, 31.11\%, 
    and 36.31\%, respectively.}
    \label{fig:fourier_recon}
\end{figure}

\Cref{fig:fourier_recon} shows the reconstructions of a single sample from the test data. Although all three reconstructions capture the overall behavior of the function, the D-MAP and greedy D-MAP estimates are closer to the truth than the Q-DEIM estimate. 
The mean relative errors of D-MAP and greedy D-MAP
over all test samples (60.66\% and 63.74\%, respectively) are close 
despite their different sensor locations.  
The mean relative error of Q-DEIM is slightly higher at 69.69\%.

In \Cref{fig:Fourier_errors}, we show the relative error of D-MAP 
and greedy D-MAP with respect to the number of
sensors $k$ and the number of POD modes $n$. 
In both panels (a) and (b), the greedy D-MAP reconstructions
maintain a relative error close to that of the D-MAP reconstructions.
This implies that the greedy D-optimal sensor locations are of similar quality
to the D-optimal sensor locations.
In panel (a), the greedy D-MAP reconstructions exhibit a slight spike in relative 
error at $n=6$.  The D-MAP reconstructions do not exhibit this spike.

The Q-DEIM relative error attains its minimum when the number of sensors and modes are equal at $n=k=5$. Clark et al. \cite{clark2021} also reported this behavior on a different data set. 
As the number of modes $n$ increases, (greedy) D-MAP errors plateau around 60\% whereas the Q-DEIM error increases.
We also note that the Q-DEIM error is always greater than or equal to that of (greedy) D-MAP.

\begin{figure}[!ht]
    \centering
    \includegraphics[width=\textwidth]{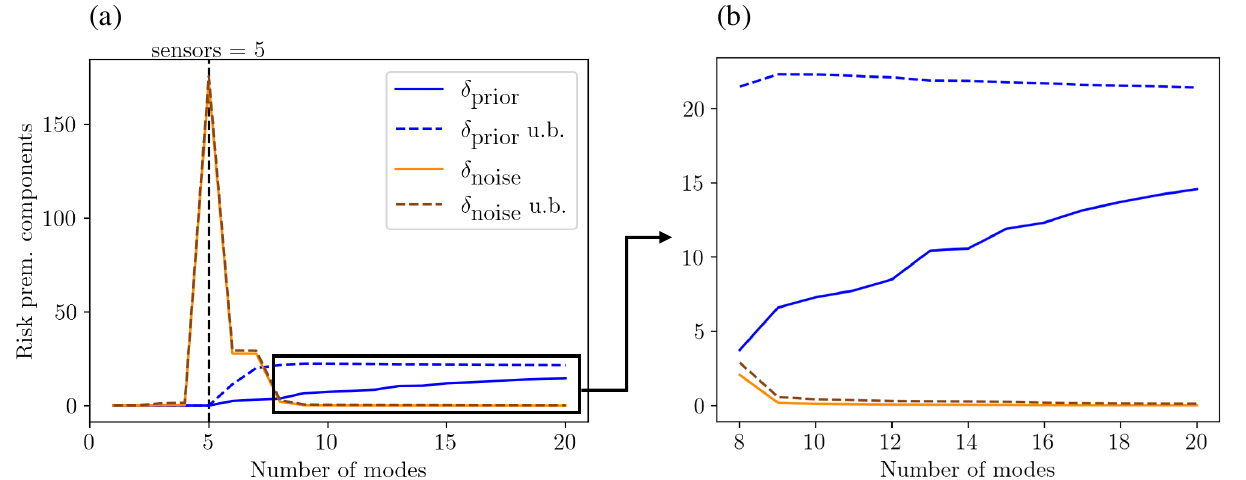}
    \caption{The components of the risk premium (i.e. $\delta_{\prior}$ and $\delta_{\noise}$) for the harmonic data and their respective upper bounds.
    We set the number of randomly-selected sensors to  $k=5$.}
    \label{fig:Fourier_risk_prem}
\end{figure}

We now turn our attention to the risk premium $\Risk(m_{\LS})-\Risk(m_{\MAP})$ to compare the reliability of Q-DEIM and MAP estimates.
In \Cref{fig:Fourier_risk_prem}, we plot the $\delta_{\prior}$ and $\delta_{\noise}$
components of the risk premium, in addition to their upper bounds, 
with respect to the number of modes $n$.
We compute this risk premium with $k=5$ sensors.  
To illustrate how an ill-conditioned $S^\top  \Phi$ affects 
the risk premium, we select these sensors 
randomly.  We vary $n$ by varying the number of
columns of the basis matrix $\Phi$ and the dimension of the prior covariance 
matrix $\Gamma_{\prior}$.
Note that $\delta_{\prior}$ equals zero when $n$ is less than
or equal to the number of sensors $k$, which
corresponds to the nullity of $S^\top  \Phi$ being 0.  When $n > k$ and
the number of modes increases, $\delta_{\prior}$ increases,
corresponding to the increased nullity of $S^\top  \Phi$.
We further note that $\delta_{\noise}$ and its upper bound $\zeta_{\noise}$ are nearly 
identical, indicating that this upper bound is tight.

The $\delta_{\noise}$ term spikes precisely when
the number of modes $n$ equals the number of sensors $k$. To explain this observation, 
let $\sigma_{\min}(S^\top  \Phi)$ denote the smallest non-zero singular 
value of $S^\top  \Phi$.  We observe that $\sigma_{\min}(S^\top  \Phi)$ drops
close to zero precisely when $S^\top \Phi$ is square, and becomes 
larger as $n$ and $k$ diverge.  A similar drop, albeit smaller in magnitude, 
occurs in the context of certain random matrices: Rudelson and Vershynin
\cite{rudelson2009randomrectangle} show that, in the case of a matrix 
$X\in \R^{k\times n}, k\geq n$ whose 
entries are iid subgaussian, the order of the smallest singular value of $X$ is at least
$\sqrt{k}-\sqrt{n-1}$ with a high probability; in the case where $X$ is square
(i.e. $X\in \R^{n\times n}$) and has iid Gaussian entries, Rudelson and Vershynin
\cite{rudelson2008randomsquare} show that, under certain moment assumptions, 
the order of the smallest singular value of $X$ is $n^{-1/2}$.  

As \cref{eq:zeta_noise_iid} indicates, a small singular value 
$\sigma_{\min}(S^\top  \Phi)$ results in a large value of $\zeta_{\noise}$.
Thus, an ill-conditioned $S^\top  \Phi$ negatively impacts the statistical quality of the 
DEIM estimate significantly more than it does the MAP estimate.
When the absolute difference between $n$ and $k$ is large, 
$\sigma_{\min}(S^\top  \Phi)$ is further from zero and $\zeta_{\noise}$ is smaller.  
Such behavior relates to the noise amplification of the DEIM
estimate described by Peherstorfer et al. \cite{peherstorfer2020}.

The upper bound of $\delta_{\prior}$ initially increases as the number of modes $n$ becomes larger
than the number of sensors $k$, then exhibits a slight decrease when $n>9$.  The upper bound of the 
risk premium exhibits a similar 
plateau. This behavior is due to the decreasing singular values
of the training data, which corresponds to decreasing eigenvalues of 
the prior covariance matrix $\Gamma_{\prior}$.
In \Cref{fig:Fourier_risk_prem}, $\delta_{\prior}$ approaches its upper bound 
as $n$ increases, corresponding to a larger portion of the variance represented by 
$\Gamma_{\prior} - \Gamma_{\post}$ being located in the null space of $S^\top  \Phi$.

\subsubsection{Q-MAP and greedy D-optimal sensors}\label{subsubsec:fourier_sensor_loc}
\begin{figure}[!ht]
    \centering
    \includegraphics[width=\textwidth]{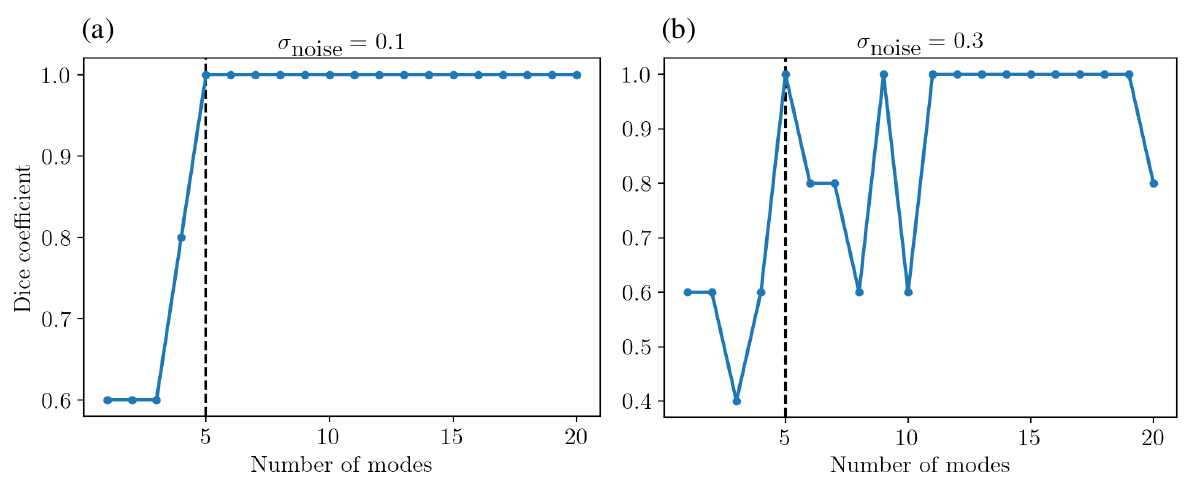}
    \caption{The Dice coefficient of the Q-MAP and greedy D-optimal sensor locations
    on the harmonic data.  We set the number of sensors $k$ to 5 (indicated by the black 
    dashed line).}
    \label{fig:Fourier_dice}
\end{figure}

In \Cref{fig:Fourier_dice}, we plot the Dice coefficient of the 
sensors selected by the greedy D-optimal algorithm and those selected by Q-MAP.  
Given a set of sensors $\xi_{\text{Q-MAP}}$ selected by Q-MAP and a set of sensors
$\xi_{\text{greedy}}$ selected by the greedy D-optimal algorithm, the
Dice coefficient \cite{dice1945} is  
\begin{equation}
    \frac{2 \, \abs{\xi_{\text{Q-MAP}} \, \cap \, \xi_{\text{greedy}}}}
    {\abs{\xi_{\text{Q-MAP}}} + \abs{\xi_{\text{greedy}}}},
\end{equation}
where $|\cdot|$ denotes the cardinality of the set.
When the sets of sensor locations selected by the two algorithms are the same, the Dice coefficient
is 1, and when the sets share no sensor locations in common, the coefficient is 0.
\Cref{fig:Fourier_dice} shows that, when noise is small and the number of modes is greater than or equal to 
the number of sensors ($n\geq k$), the sensor locations selected by the Q-MAP and 
greedy D-optimal algorithms are identical.

Note the lack of correspondence between the Q-MAP and greedy D-optimal sensor locations 
when the number of sensors $k$ exceeds the number of modes $n$.  This stems from the 
CPQR algorithm only being able to place $n$ sensors, with the remaining $k-n$ sensors
placed arbitrarily. 

Now consider the case when $k\leq n$.
Recall that, as a scalar multiple of the identity, $\Gamma_{\noise}$ has no impact on 
the Q-MAP sensor locations.  
Such is not the case for the greedy D-optimal sensor locations.  When the 
effect of noise regularization is significant, the Q-MAP and greedy D-optimal sensor 
locations can be distinct, as we see with the dips in \Cref{fig:Fourier_dice} (b) 
that are not apparent in \Cref{fig:Fourier_dice} (a).  
Nonetheless, the close correspondence between the sensor locations of these two 
algorithms, particularly in the case of a small $\sigma_{\noise}$, is remarkable.

\subsection{Two-dimensional turbulence}\label{subsec:turb}
Our second example concerns a two-dimensional fluid flow whose vorticity field $\omega(x,y,t)$ satisfies
\begin{equation}\label{eq:turb}
    \partial_t \omega + \mathbf{u} \cdot \nabla \omega 
    = \nu\, \Delta \omega - n_f \cos (n_f\, \mathbf{y}).
\end{equation}
Denoting the stream function by $\phi(x,y,t)$, the velocity field is given by
$\mathbf{u}=(\partial_{y} \psi, -\partial_{x} \psi)$, and the vorticity field satisfies $\omega = -\Delta \psi$.  
The integer $n_f=4$ denotes the forcing wavenumber and $\nu=1/40$
is the inverse of the Reynolds number.  The spatial domain is 
$(x,y)\in[0, 2\pi] \times [0, 2\pi]$ with periodic boundary conditions.  
The resulting flow is referred to as Kolmogorov flow and is known to exhibit chaotic dynamics~\cite{faraz_adjoint,Farazmande1701533}.

We solve equation \Cref{eq:turb} using a standard pseudo-spectral numerical  method, over an equispaced spatial grid of size $128 \times 128$.
As a result, the high-fidelity resolution is $N=128^2$.  For additional details regarding this implementation,
we refer to Section 3.1 of Farazmand and Saibaba~\cite{farazmand2023}.

We gather 1001 snapshots in time of the vorticity field $\omega(x,y,t_i)$, $i=1,2,\cdots,1001$.
Snapshots are only taken after sufficient time has elapsed to remove the 
effects of the initial transients.  We retain the first 750 snapshots of $\omega$ for training, 
and the remaining 251 snapshots for testing.
The first 100 POD modes account for approximately 95\% of training data variance, which is why we set $n=100$ in our results unless explicitly stated otherwise.
We add iid Gaussian noise to  the test data with mean $0$ and standard deviation $\sigma_{\noise}=0.3$, so that
$\Gamma_{\noise}=\sigma_{\noise}^2 I$.
On average, noise represents approximately 15\% of the test data.

\subsubsection{Reconstruction error and the risk premium}\label{subsubsec:turb_recon}
\begin{figure}[!ht]
    \centering
    \includegraphics[width=.85\textwidth]{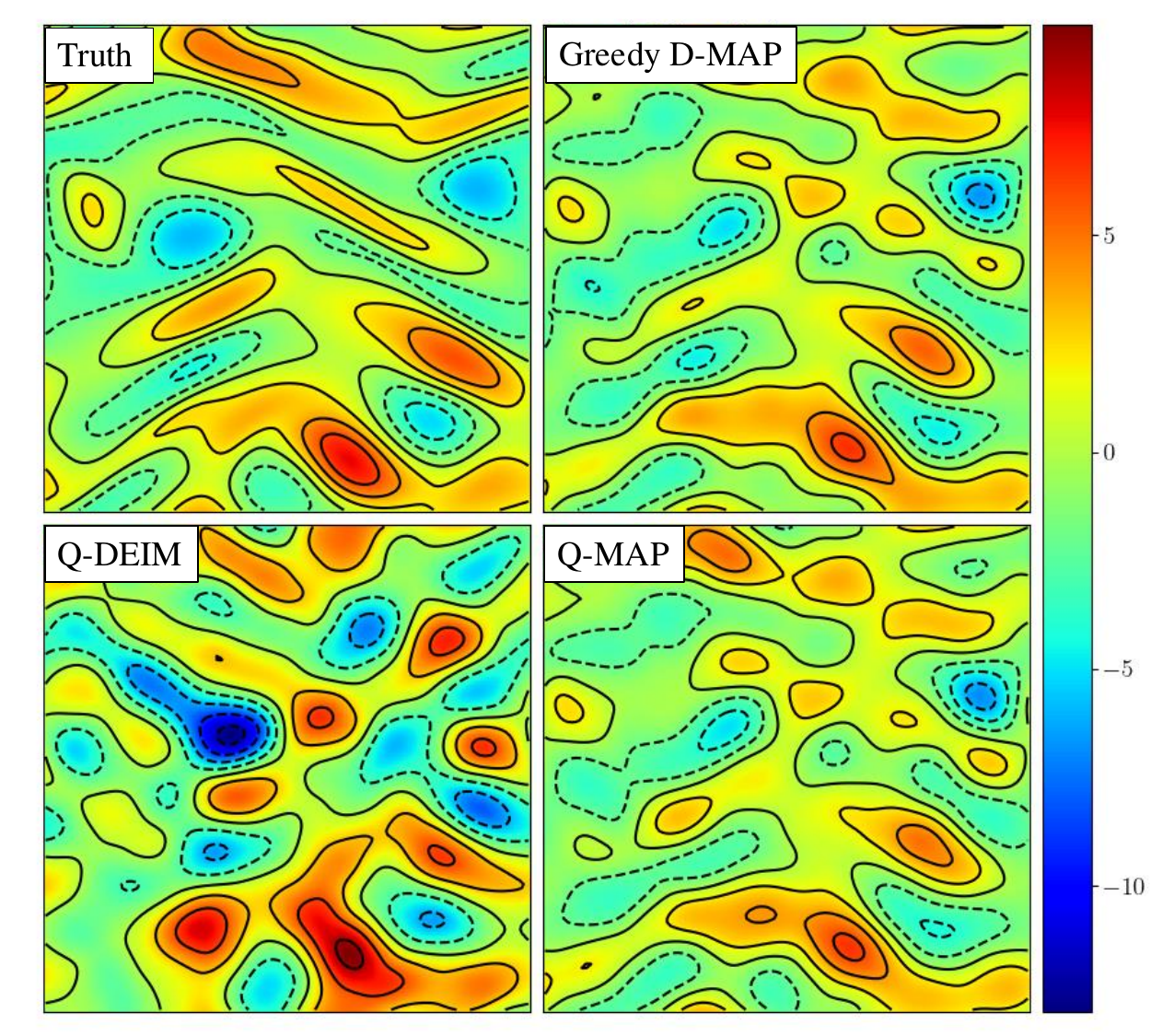}
    \caption{Plot of the reconstructions of a single turbulence data sample.
    We set the number of sensors $k$ to 50.  For Q-DEIM, we set the number of
    modes $n$ equal to $k$ (50), while for greedy D-MAP
    and Q-MAP we set $n=100$.  On this sample, greedy D-MAP,
    Q-DEIM, and Q-MAP have relative errors of 54.05\%, 107.67\%, 
    and 50.88\%, respectively.}
    \label{fig:turb_recon}
\end{figure}

\Cref{fig:turb_recon} shows the reconstructions of a snapshot from the test
samples using 50 sensors. As with the harmonic data, the greedy D-MAP reconstruction is more accurate than
the Q-DEIM reconstruction. More precisely, the relative error of greedy D-MAP is about 54\%, whereas the relative error of Q-DEIM is 108\%.
The Q-MAP reconstruction is distinct from the greedy D-MAP reconstruction, but its relative error is similar (approximately 51\%).

\begin{figure}[!ht]
    \centering
    \vspace*{-2.5cm}
    \hspace*{-1cm}
    \includegraphics[width=\textwidth]{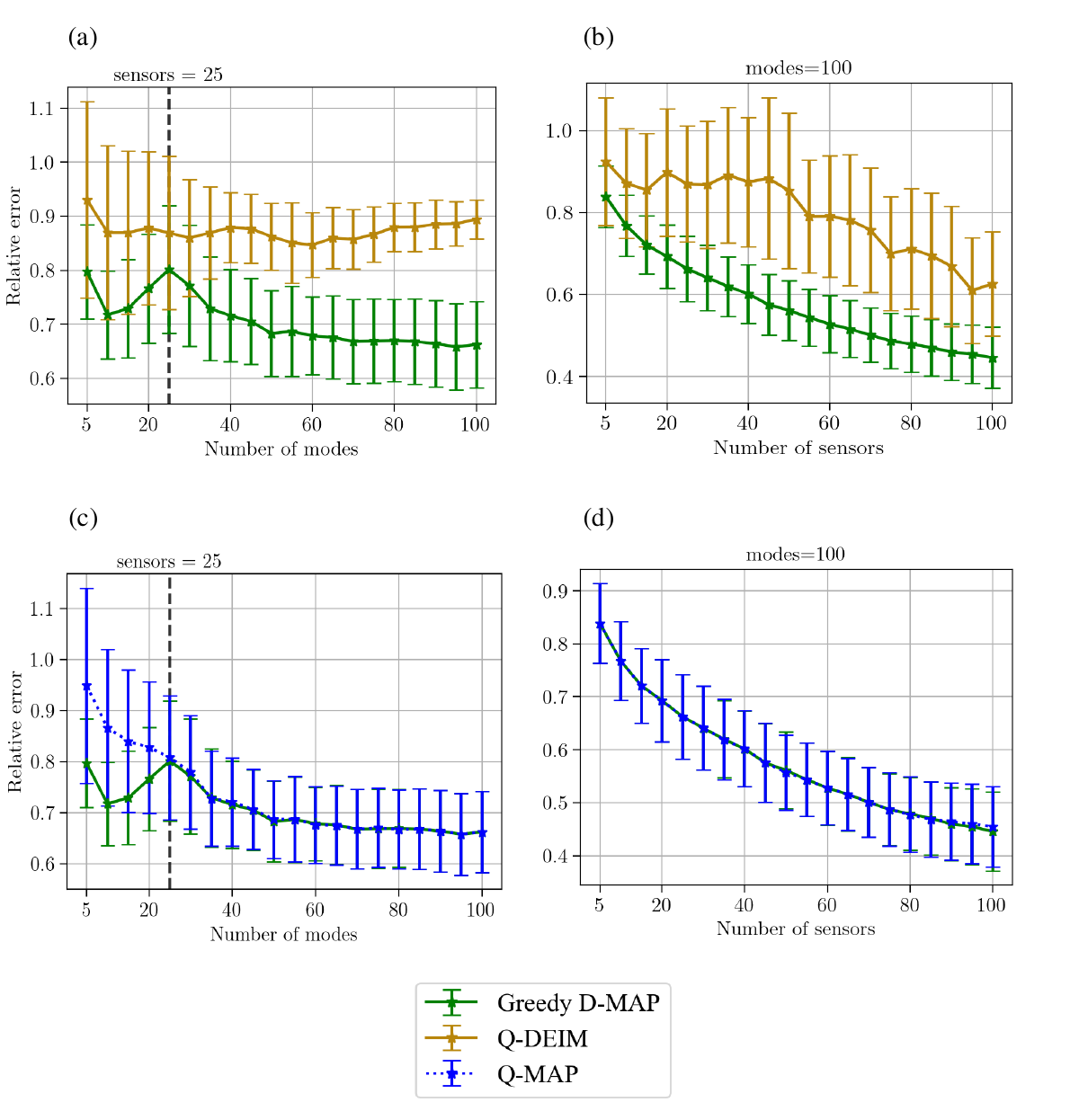}
    \caption{We plot the relative error of greedy D-MAP and Q-DEIM
    on the turbulence data with respect to number of modes in (a) and sensors in (b).  We plot the
    relative error of greedy D-MAP and Q-MAP with respect to number of modes
    in (c) and sensors in (d).  In (b), we set the number of modes used by Q-DEIM to the number
    of sensors (on the horizontal axis), and the number of modes used by greedy D-MAP to 100.}
    \label{fig:turb_errors}
\end{figure}

Now we turn our attention to examining the relative errors over the entire test dataset. In \Cref{fig:turb_errors}, we see greedy D-MAP 
outperforming Q-DEIM over a wide range of number of modes $n$ and sensors $k$.  In panel (a), as the number of modes $n$ increases past the number of
sensors $k$, the mean Q-DEIM error initially exhibits no significant upward or downward trend, 
then exhibits a slight upward trend when $n>60$.  
In contrast, the greedy D-MAP error decreases when $25 < n \leq 70$, 
and plateaus when $n> 70$.  As $n$ approaches $k$, the greedy D-MAP error spikes.  
The spike in this case is more distinct than that of \Cref{fig:Fourier_errors} (a).
In panel (b), we see that
as the number of sensors $k$ increases, the Q-DEIM error initially exhibits no significant upward or downward
trend, then decreases when $k> 45$.  In contrast, the
greedy D-MAP error decreases monotonically.  Panels (c) and (d) show that,
when $n\geq k$, greedy D-MAP and Q-MAP have identical or nearly identical relative errors. However, when the number of modes is smaller than the number of sensors ($n<k$), greedy D-MAP outperforms Q-MAP.

\begin{figure}[!ht]
    \centering
    \includegraphics[width=\textwidth]{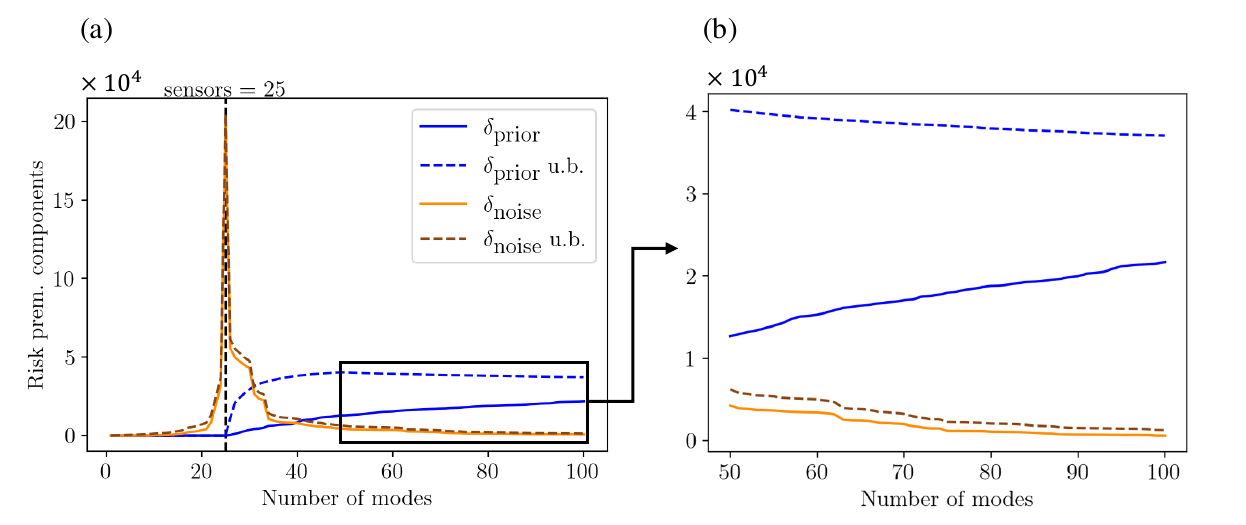}
    \caption{(a) and (b) plot the components of the turbulence data risk premium,
    in addition to their respective upper bounds.  
    We set the number of randomly-selected sensors $k$ to 25.}
    \label{fig:turb_risk_prem}
\end{figure}

\Cref{fig:turb_risk_prem} shows the components of the risk premium and their 
upper bounds as a function of the number of modes $n$.  We compute this risk premium from 25 randomly selected 
sensors.  The shapes of the risk premium components and their upper bounds are similar to 
\Cref{fig:Fourier_risk_prem} for the harmonic data.

In both cases (\Cref{fig:Fourier_risk_prem,fig:turb_risk_prem}), there is a spike in the $\delta_{\noise}$ term
which occurs when the number of modes $n$ 
equals the number of sensors $k$.  Since the risk premium equals 
$\delta_{\prior}+\delta_{\noise}$, a spike in the $\delta_{\noise}$ term corresponds to a 
large risk premium.  Thus, the Bayes risk of the MAP estimate is significantly lower 
than the Bayes risk of the MNLS estimate when $n=k$, suggesting that the MAP estimate performs 
significantly better.  However, in Figures \ref{fig:Fourier_errors} (c) and \ref{fig:turb_errors} (a), which show the relative error,
there is no significant difference between the greedy D-MAP and Q-DEIM error
when $n=k$.

This seeming discrepancy arises from the different sensor placement methods used. In our plots of the risk premium, we place $k$ sensors randomly at fixed positions.
These sensor locations do not change as the number of modes $n$ increases. On the other hand, in the figures showing the relative error, the optimal sensor locations are estimated and updated for each value of $n$. For instance, for Q-DEIM, the CPQR algorithm is used to compute the sensor locations.
When the number of modes equals the number of sensors ($n=k$), Drma\v{c} and Gugercin \cite{drmac2016} 
proved that using CPQR to select sensors results in an upper bound on $\norm{(S^\top  \Phi)^{-1}}_2=1/\sigma_{\min}(S^\top  \Phi)$ 
or, equivalently, a lower bound on $\sigma_{\min}(S^\top  \Phi)$.  This lower bound ensures that $\sigma_{\min}(S^\top  \Phi)$ 
is not too close to zero and consequently that $\delta_{\noise}$ is not too large.
Thus, while the risk premium spikes when sensors are chosen randomly,  using CPQR to place sensors mitigates this spike.

\subsubsection{Q-MAP and greedy D-optimal sensors}\label{subsubsec:turb_sensor_loc}
\begin{figure}[!ht]
    \centering
    \includegraphics[width=\textwidth]{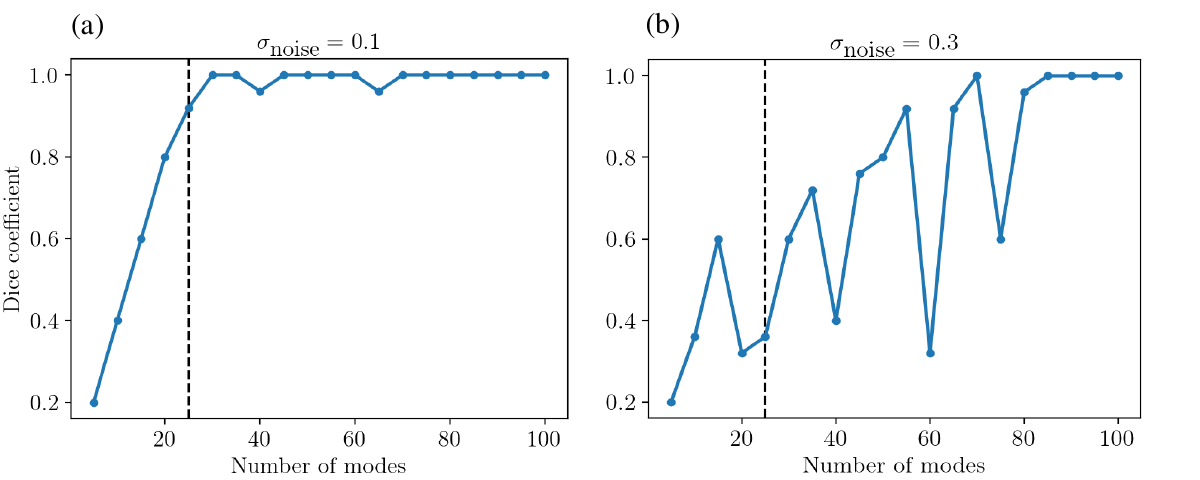}
    \caption{The Dice coefficient of the Q-MAP and greedy D-optimal sensor locations
    on the turbulence data.  We set the number of sensors $k$ to 25 (indicated by the black 
    dashed line).}
    \label{fig:turb_dice}
\end{figure}

\Cref{fig:turb_dice} shows the Dice coefficient for the turbulence data.
Similar to the harmonic functions (\Cref{fig:Fourier_dice}), when measurement noise is sufficiently small and the 
number of modes $n$ is greater than or equal to the number of sensors $k$, the greedy
D-optimal sensor locations and those selected by Q-MAP are identical or nearly identical.
For larger measurement noise, \Cref{fig:turb_dice} (b) shows that the greedy D-optimal and CPQR algorithms return different sensor locations. 
Despite these differences, 
\Cref{fig:turb_errors} (c-d) show that Q-MAP and greedy D-MAP have nearly identical relative errors when $n \geq k$.

\section{Conclusion}\label{sec:conc}
Estimating the full state of a system from partial observations has applications in many areas of science such as meteorology, control theory, and fluid dynamics.
We considered two methods that seek to accomplish this task: a probabilistic Bayesian method (the MAP estimate) and a deterministic least-squares method (DEIM). 

Theoretical results involving the relative error of these methods are not informative since computable expressions for the errors are unavailable. As we discussed in \Cref{subsubsec:map_rel_err}, upper bounds on these errors also fail to provide illuminating information regarding their relative accuracy. As a result, we proposed the Bayesian risk premium as
an alternative criterion to compare the quality of the MAP and DEIM estimates. We derived a computable expression for the risk premium, which indicates that the MAP estimate is always more reliable than the DEIM estimate. Furthermore, we proved that the risk premium has two nonnegative components. One component quantifies prior uncertainty ($\delta_{\prior}$), whereas the other reflects the effects of measurement noise ($\delta_{\noise}$). 

Our numerical results indicate that, when the number of modes exceeds the number of sensors, $\delta_{\prior}$ dominates while $\delta_{\noise}$
remains close to zero. In contrast, when the number of modes equals the number of sensors, $\delta_{\noise}$ spikes due to ill-conditioning of the observation matrix $S^\top \Phi$.
This helps explain the large relative reconstruction error observed by Klishin et al. \cite{klishin2023} and Clark et al. \cite{clark2021} in this regime.

Although informative theoretical results regarding the relative reconstruction errors are not available, our numerical results returned smaller relative errors for the MAP estimate compared to the DEIM estimate. These numerical observations are consistent with our theoretical results based on the risk premium.

Accuracy of MAP estimates depend to a great extent on the choice of the prior distribution. We proposed a new prior covariance matrix which is derived systematically 
from a change of basis applied to the sample covariance matrix. This prior produces reconstructions whose relative errors decrease monotonically as the number of sensors increases.

Looking forward to the future, our theoretical results motivate new optimal sensor placement methods.
Remarkably, the only component of the prior covariance matrix relevant to the Bayes risk 
of the MNLS estimator is that which lies in the null space of $S^\top \Phi$.
This is the variance of the inversion parameter $m$ that does not result in variance in the 
measurements $y$.  We showed that the magnitude of the posterior 
covariance in the null space of $S^\top \Phi$ is never greater than the magnitude
of the prior covariance in the null space of $S^\top \Phi$, meaning variance
in this subspace never increases after taking measurements.
With this result in mind, future sensor placement algorithms could 
seek to minimize,
\begin{equation}
\tr[\Gamma_{\prior}-(S^\top\Phi)^+ (S^\top\Phi) \Gamma_{\prior}],
\label{eq:tr_proj}
\end{equation}
which represents the magnitude of the columns of $\Gamma_{\prior}$ orthogonally projected onto the 
null space of $S^\top \Phi$.
Recent results~\cite{farazmand2024, farazmand2025} have shown that the component of the inversion 
parameter $m$ that lies in the null space of 
$S^\top \Phi$ has a significant impact on reconstruction quality. 
Sensor placements that minimize~\eqref{eq:tr_proj} would be minimizing the variance 
of this component of the inversion parameter.
When $S^\top \Phi$ is nearly square, the impact of measurement noise on the quality of the DEIM 
estimate is amplified.  In this case, a natural choice of 
design criterion is \Cref{eq:ls_bayes_risk}, which is the Bayes risk of the MNLS 
estimator and an upper bound on the Bayes risk of the MAP estimator.

\vspace{1em}
\noindent\textbf{Funding.} {This work was supported by 
the National Science Foundation (NSF) under grant DMS-1745654. 
AA was supported in part by NSF award DMS-2111044.
MF was supported in part by NSF award DMS-2220548.
AKS was supported in part by NSF award DMS-1845406 and
the Department of Energy through the award DE-SC0023188. 
}

\vspace{1em}
\noindent\textbf{Data availability statement.} Code and data are available in the following GitHub repository.
{\\\href{https://github.com/LevKakasenko/state-estimation-bridge}
{https://github.com/LevKakasenko/state-estimation-bridge}}

\appendix
\renewcommand{\thesection}{\Alph{section}}
\makeatletter
\renewcommand\@seccntformat[1]{\appendixname\ \csname the#1\endcsname.\hspace{0.5em}}
\makeatother

\section{Proof of upper bounds on the MAP and DEIM relative errors}\label{appendix:1}
\begin{proof}[Proof of \Cref{prop:map_ub}]
	Note that the relative error between the true full state $u$
	and the MAP reconstruction $u_{\MAP}$ is
	\begin{align*}
		\begin{split}\label{eq:map_ub_}
			\frac{\norm{u_{\MAP}-u}_2}{\norm{u}_2} 
			=&\frac{\norm{\Phi m_{\MAP}-u}_2}{\norm{u}_2} \\
			=&\frac{\norm{\Phi \Gamma_{\post}(S^\top  \Phi)^\top 
					\Gamma_{\noise}^{-1}y-u}_2}{\norm{u}_2} \\
			=&\frac{\norm{\Phi \Gamma_{\post}(S^\top  \Phi)^\top 
					\Gamma_{\noise}^{-1}(S^\top u+\eta)-u}_2}{\norm{u}_2} \\
			\leq& \norm{\Phi \Gamma_{\post}(S^\top  \Phi)^\top 
				\Gamma_{\noise}^{-1}S^\top -I}_2 \\
			&+\norm{\Gamma_{\post}(S^\top  \Phi)^\top 
				\Gamma_{\noise}^{-1}}_2\, \frac{\norm{\eta}_2}{\norm{u}_2}
			\period
		\end{split}
	\end{align*}
\end{proof}

\begin{proof}[Proof of \Cref{prop:deim_ub}]
	Note that the relative error between the true full state $u$
	and the DEIM reconstruction $u_{\DEIM}$ is
	\begin{align*}
		\begin{split}
			\frac{\norm{u_{\DEIM}-u}_2}{\norm{u}_2} 
			=&\frac{\norm{\Phi m_{\LS}-u}_2}{\norm{u}_2} \\
			=& \frac{\norm{\Phi(S^\top  \Phi)^+ y - u}_2}{\norm{u}_2} \\
			=& \frac{\norm{\Phi(S^\top  \Phi)^+ 
					(S^\top u + \eta) - u}_2}{\norm{u}_2} \\
			\leq& \norm{\Phi(S^\top  \Phi)^+S^\top -I}_2 \\
			&+\norm{(S^\top  \Phi)^+}_2\, \frac{\norm{\eta}_2}{\norm{u}_2}
			\period
		\end{split}
	\end{align*}
	Since $P:=\Phi(S^\top  \Phi)^+S^\top $ is a projection that is not null
	and not the identity, we have that 
	$\norm{P-I}_2=\norm{P}_2$ \cite{chatur2010, szyld2006}.  Thus,
	\begin{align*}
		\begin{split}
			\frac{\norm{u_{\DEIM}-u}_2}{\norm{u}_2} 
			\leq& \norm{P}_2
			+\norm{(S^\top  \Phi)^+}_2\, \frac{\norm{\eta}_2}{\norm{u}_2} \\
			=& \norm{(S^\top  \Phi)^+}_2
			+\norm{(S^\top  \Phi)^+}_2\, \frac{\norm{\eta}_2}{\norm{u}_2}
			\period
		\end{split}
	\end{align*}
\end{proof}

\section{Proof of \Cref{thm:bayes_risk_ineq}}\label{appendix:2}
The following are two supporting lemmas, followed by the proof of \Cref{thm:bayes_risk_ineq}.
\begin{lemma}\label{lem:spsd1}
	Let $L\in \R^{n \times n}$ be some symmetric positive semidefinite matrix.  
	Then $I-(I+L)^{-1}$ is symmetric positive semidefinite.
\end{lemma}
\begin{proof}
	Let $H:=I-(I+L)^{-1}$, and $L$ have the eigen-decomposition 
	$L=V\Lambda V^{-1}$, where 
	$\Lambda=\diag(\lambda_1, ..., \lambda_n)$.  Then
	\begin{align*}
		\begin{split}
			H &= I-(I+L)^{-1} \\
			&= I-(I+V\Lambda V^{-1})^{-1} \\
			&= V \underbrace{\sqbr{I-(I+\Lambda)^{-1}}}_{:=M}V^{-1} \text{ .}
		\end{split}
	\end{align*}
	
	Note that $M$ is a diagonal matrix whose entries are the eigenvalues of $H$.  
	The $i^{\th_}$ diagonal entry of $M$ is
	\begin{equation*}
		m_{i,i}=1-\frac{1}{1+\lambda_i}\geq 0
	\end{equation*}
	since $\lambda_i \geq 0$. 
	
	Because $H$ is symmetric with non-negative eigenvalues, $H$ is symmetric positive semidefinite.
\end{proof}

\begin{lemma}\label{lem:spsd2}
	Let $B\in \R^{k \times n}$ and $P:=I-B(I+B^\top B)^{-1}B^\top $.  Then $P$ is symmetric positive definite.
\end{lemma}
\begin{proof}
	Let $B$ have the full singular value decomposition $B=U\Sigma V^\top $.  Then
	\begin{align*}
		\begin{split}
			P&=I-B(I+B^\top B)^{-1}B^\top  \\
			&=U\underbrace{\sqbr{I-\Sigma (I+\Sigma^\top \Sigma)^{-1} \Sigma^\top }}_{:=W} U^\top  \text{ .}
		\end{split}
	\end{align*}
	Note that $W$ is a diagonal matrix whose entries are the eigenvalues of $P$.  
	If $B$ has rank $r$, then the first $r$ diagonal entries of $W$ are
	\begin{equation*}
		w_{i,i}=1-\frac{\sigma_i^2}{1+\sigma_i^2} > 0 \text{ ,}
	\end{equation*}
	where $\sigma_i$ is the $i^{\th_}$ singular value of $B$.  
	The remaining $k-r$ diagonal entries of $W$ are 1.  
	
	Since $P$ is symmetric with strictly positive eigenvalues, $P$ is symmetric positive definite.
\end{proof}

We are now ready to the state the proof of \Cref{thm:bayes_risk_ineq}.

\begin{proof}[Proof of \Cref{thm:bayes_risk_ineq}]
	Note that
	\begin{align}
		\begin{split}
			\Risk(m_{\MAP}):&=\E_m[E_{y|m}[\norm{m_{\MAP}-m}_2^2]] \\
			&=\tr[(\Gamma_{\prior}^{-1}+A^\top \Gamma_{\noise}^{-1}A)^{-1}] \text{ ,} \\
			\Risk(m_{\LS}):&=\E_m[E_{y|m}[\norm{m_{\LS}-m}_2^2]] \\
			&= \tr[(I-A^+A)\Gamma_{\prior}+A^+\Gamma_{\noise}(A^+)^\top ] \text{ .}
		\end{split}
	\end{align}
	We can re-express $\Risk(m_{\MAP}), \Risk(m_{\LS})$ as
	\begin{align}\label{eq:r_ls__r_map}
		\begin{split}
			\Risk(m_{\MAP})=&\underbrace{\tr[(I-A^+A)(\Gamma_{\prior}^{-1}+A^\top \Gamma_{\noise}^{-1}A)^{-1}]}_{:=\alpha_{\prior}} \\
			&+\underbrace{\tr[A^+A(\Gamma_{\prior}^{-1}+A^\top \Gamma_{\noise}^{-1}A)^{-1}]}_{:=\alpha_{\noise}} \text{ ,} \\
			\Risk(m_{\LS})=&\underbrace{\tr[(I-A^+A)\Gamma_{\prior}]}_{:=\beta_{\prior}}+\underbrace{\tr[A^+\Gamma_{\noise}(A^+)^\top ]}_{:=\beta_{\noise}} \text{ .}
		\end{split}
	\end{align}
	
	Consider $\beta_{\prior} - \alpha_{\prior}$.  Note that
	\begin{align}\label{eq:beta1_min_alpha1}
		\begin{split}
			\delta_{\prior} &= \beta_{\prior} - \alpha_{\prior} \\
			&= \tr \sqbr{(I-A^+A)\left(\Gamma_{\prior}-(\Gamma_{\prior}^{-1}+A^\top \Gamma_{\noise}^{-1}A)^{-1}\right)} \\
			&= \tr \sqbr{(I-A^+A)\Gamma_{\prior}^{1/2}
				(I-(I+B^\top B)^{-1})\Gamma_{\prior}^{1/2}} \com
		\end{split}
	\end{align}
	where $B:=\Gamma_{\noise}^{-1/2}A\Gamma_{\prior}^{1/2}$.  Note 
	that $B^\top B$ is symmetric positive semi-definite (spsd).  
	Thus, by \Cref{lem:spsd1}, we have that $I-(I+B^\top B)^{-1}$ is spsd.  
	Thus, $I-(I+B^\top B)^{-1}=K^\top K$ for some matrix $K$.  
	
	Also note that $I-A^+A$ is an orthogonal projection, and thus spsd.  We then have that $I-A^+A=F^\top F$ for some matrix $F$.
	
	Substituting into \Cref{eq:beta1_min_alpha1}, we get
	\begin{align*}
		\begin{split}
			\delta_{\prior} &= \tr \left(F^\top F \Gamma_{\prior}^{1/2}K^\top K \Gamma_{\prior}^{1/2} \right) \\
			&= \tr \sqbr{(K\Gamma_{\prior}^{1/2}F^\top )^\top (K\Gamma_{\prior}^{1/2}F^\top )} \\
			&\geq 0 \text{ .}
		\end{split}
	\end{align*}
	
	Now consider $\beta_{\noise} - \alpha_{\noise}$.  Note that
	\begin{align}\label{eq:beta2_min_alpha2}
		\begin{split}
			\delta_{\noise} &= \beta_{\noise} - \alpha_{\noise} \\
			&= \tr \sqbr{ A^+ \Gamma_{\noise} (A^+)^\top  - A^+A (\Gamma_{\prior}^{-1}+A^\top \Gamma_{\noise}^{-1}A)^{-1}} \\
			&= \tr \sqbr{(A^+ \Gamma_{\noise}^{1/2}) \underbrace{\sqbr{I-B(I+B^\top B)^{-1}B^\top }}_{:=P} (A^+ \Gamma_{\noise}^{1/2})^\top }
			\period
		\end{split}
	\end{align}
	By \Cref{lem:spsd2}, we have that matrix $P$ is symmetric positive definite.  
	Thus, $P=Q^\top Q$ for some matrix $Q$.
	Substituting into \Cref{eq:beta2_min_alpha2}, we get
	\begin{align*}
		\begin{split}
			\delta_{\noise} &= \tr \sqbr{(A^+ \Gamma_{\noise}^{1/2})Q^\top Q(A^+ \Gamma_{\noise}^{1/2})^\top } \\
			&= \tr \sqbr{(A^+ \Gamma_{\noise}^{1/2}Q^\top )^\top (A^+ \Gamma_{\noise}^{1/2}Q^\top )} \\
			&\geq 0 \text{ .}
		\end{split}
	\end{align*}
	
	By \Cref{eq:r_ls__r_map} and the definitions of $\delta_{\prior}$ and $\delta_{\noise}$, 
	we have
	\begin{equation}
		\Risk(m_{\LS}) - \Risk(m_{\MAP})=\delta_{\prior}+\delta_{\noise} \period
	\end{equation}
\end{proof}

\section{Proof of \Cref{thm:risk_prem_ub}}\label{appendix:3}
The following are two supporting lemmas, followed by the proof of \Cref{thm:risk_prem_ub}.

\begin{lemma}\label{lem:delta_prior_ub}
	The $\delta_{\prior}$ component of the risk premium is bounded
	as
	\begin{equation}
		\delta_{\prior}
		\leq \sum_{i=1}^{\textnormal{Nullity}\paren{A}}
		\lambda_i(\Gamma_{\prior}-\Gamma_{\post})
		\period
	\end{equation}
	This bound on $\delta_{\prior}$ is tight.
\end{lemma}
\begin{proof}
	By the Von Neumann trace inequality,
	\begin{align}\label{eq:delta_prior_ub_}
		\begin{split}
			\delta_{\prior}&=\tr \sqbr{\paren{I-A^+ A} 
				\paren{\Gamma_{\prior}-\Gamma_{\post}}} \\
			&\leq \sum_{i=1}^{n} \lambda_i \paren{I-A^+ A}
			\lambda_i \paren{\Gamma_{\prior}-\Gamma_{\post}} \period
		\end{split}
	\end{align}
	Since $I-A^+ A$ is an orthogonal projection onto the null space of
	$A$, its eigenvalues are either 0 or 1, and its rank is the nullity
	of $A$.  Thus,
	\begin{equation}
		\delta_{\prior}
		\leq \sum_{i=1}^{\textnormal{Nullity}\paren{A}}
		\lambda_i(\Gamma_{\prior}-\Gamma_{\post}) \period
	\end{equation}
	
	\paragraph{We now show that this bound on $\delta_{\prior}$ is tight}
	Note that for any real symmetric matrices $A$ and $B$ of the same
	dimension, 
	\begin{equation}
		\tr \paren{AB}=\sum_{i=1}^{n} \lambda_i(A)\lambda_i(B)
	\end{equation}
	if and only if $A$ and $B$ have the respective eigen-decompositions
	\begin{align}
		A&=V\Lambda_A V^\top  \com \\
		B&=V\Lambda_B V^\top  \com
	\end{align}
	where $\Lambda_A$ and $\Lambda_B$ are diagonal matrices sorted in
	descending order \cite{carlsson2021}.  Thus, equality holds in 
	\Cref{eq:delta_prior_ub_} if and only if the eigenvectors
	corresponding to the largest eigenvalues of $I-A^+A$ are also the
	eigenvectors corresponding to the largest eigenvalues of 
	$\Gamma_{\prior}-\Gamma_{\post}$.  The largest possible eigenvalues 
	of $I-A^+A$ are equal to 1, with corresponding eigenvectors
	in the null space of $A$.
	
	Thus, equality holds in \Cref{eq:delta_prior_ub_} if
	and only if the eigenvectors corresponding to the largest 
	eigenvalues of $\Gamma_{\prior}-\Gamma_{\post}$ are all in the 
	null space of $A$.
\end{proof}

\begin{lemma}\label{lem:delta_noise_ub}
	The  $\delta_{\noise}$ component of the risk premium is bounded
	above by
	\begin{equation}
		\delta_{\noise} \leq \tr \sqbr{A^+ \Gamma_{\noise} (A^+)^\top } \period
	\end{equation}
\end{lemma}
\begin{proof}
	Note that
	\begin{equation}
		\delta_{\noise}=\tr \sqbr{A^+ \Gamma_{\noise} (A^+)^\top } - \tr \paren{A^+A \Gamma_{\post}} \period
	\end{equation}
	By the Von Neumann trace inequality, we have that
	\begin{equation}
		\tr \paren{A^+A \Gamma_{\post}} \geq \sum_{i=1}^{n} \lambda_i \paren{A^+A} \lambda_{n-i+1} 
		\paren{\Gamma_{\post}} \period
	\end{equation}
	Since the eigenvalues of $A^+A$ and $\Gamma_{\post}$ are nonnegative, we have that
	\begin{equation}
		\tr \paren{A^+A \Gamma_{\post}} \geq 0 \com
	\end{equation}
	which implies that 
	\begin{align}
		\begin{split}
			\delta_{\noise}&=\tr \sqbr{A^+ \Gamma_{\noise} (A^+)^\top } - \tr \paren{A^+A \Gamma_{\post}} \\
			&\leq \tr \sqbr{A^+ \Gamma_{\noise} (A^+)^\top } \period
		\end{split}
	\end{align}
\end{proof}

\begin{proof}[Proof of \Cref{thm:risk_prem_ub}]
	By \Cref{thm:bayes_risk_ineq}, we have that
	\begin{equation}
		0 \leq \Risk(m_{\LS}) - \Risk(m_{\MAP}) 
		= \delta_{\prior} + \delta_{\noise} \period
	\end{equation}
	By \Cref{lem:delta_prior_ub}, we have that
	\begin{equation}
		\delta_{\prior}
		\leq \sum_{i=1}^{\textnormal{Nullity}\paren{A}}
		\lambda_i(\Gamma_{\prior}-\Gamma_{\post})
	\end{equation}
	By \Cref{lem:delta_noise_ub}, we have that
	\begin{equation}
		\delta_{\noise} \leq \tr \sqbr{A^+ \Gamma_{\noise} (A^+)^\top } \period
	\end{equation}
	Thus,
	\begin{align}
		\begin{split}
			0 &\leq \Risk(m_{\LS}) - \Risk(m_{\MAP}) \\
			&\leq \tr \sqbr{A^+ \Gamma_{\noise} (A^+)^\top } +
			\sum_{i=1}^{\textnormal{Nullity}\paren{A}} \lambda_i(\Gamma_{\prior}-\Gamma_{\post})
			\period
		\end{split}
	\end{align}
\end{proof}


\end{document}